\documentclass[11pt]{article}
\usepackage{amsthm,amsmath,amssymb,bbm,geometry,epsfig,hyperref,enumerate,comment,nicefrac,listings,graphicx,xcolor,color,titling,makecell}
\usepackage{datetime}
\definecolor{codeblue}{rgb}{0,0,0.8}
\definecolor{codegray}{rgb}{0.5,0.5,0.5}
\definecolor{codebrown}{rgb}{0.56,0.28,0.16}
\definecolor{backcolour}{rgb}{1,1,1}
\definecolor{codegreen}{rgb}{0.4,0.8,0.3}

\hypersetup{
    colorlinks,
    linkcolor={blue!80!black},
    citecolor={blue},
    urlcolor={blue!80!black}
}

\lstdefinestyle{mystyle}{
	language = Python,
	morekeywords = {as, with}, 
	backgroundcolor=\color{backcolour},   
	keywordstyle=\color{codeblue},
	commentstyle=\color{codegray},
	numberstyle=\tiny,
	stringstyle=\color{codegreen},
	basicstyle=\footnotesize,
	breakatwhitespace=false,         
	breaklines=true,                 
	captionpos=b,                    
	keepspaces=true,                 
	numbers=left,  
	frame = single,                  
	numbersep=5pt,                  
	showspaces=false,                
	showstringspaces=false,
	showtabs=false,                  
	tabsize=2
}
\newtheorem{lemma}{Lemma}[section]

\newtheorem{algo}[lemma]{Framework}
\theoremstyle{definition}
\newtheorem{remark}{Remark}[section]

 \newcommand{\be}{\begin{equation}}
 \newcommand{\ee}{\end{equation}}
 \newcommand{\bea}{\begin{eqnarray}}
 \newcommand{\eea}{\end{eqnarray}}
 \newcommand{\beas}{\begin{eqnarray*}}
\newcommand{\eeas}{\end{eqnarray*}}

\newcommand{\edg}[1]{\ensuremath{\! \left[ #1 \right] }}
\newcommand{\brak}[1]{\ensuremath{\left( #1 \right)}}

\newcommand{\abs}[1]{\ensuremath{ \left| #1 \right| }}

\providecommand{\N}{{\ensuremath{\mathbb{N}}}}

\providecommand{\R}{{\ensuremath{\mathbb{R}}}}
\providecommand{\B}{\mathcal{B}}

\renewcommand{\P}{\mathbbm{P}}

\providecommand{\bS}{\mathbb{S}}

\providecommand{\cY}{\mathcal{Y}}

\providecommand{\cV}{\mathcal{V}}
\renewcommand{\S}{\mathcal{S}}
\providecommand{\E}{{\ensuremath{\mathbbm{E}}}}

\providecommand{\N}{{\ensuremath{\mathbbm{N}}}}
\providecommand{\bV}{{\ensuremath{\mathbb{V}}}}

\providecommand{\R}{{\ensuremath{\mathbbm{R}}}}

\providecommand{\E}{{\ensuremath{\mathbb{E}}}}

\newcommand{\F}{{\ensuremath{\mathcal{F}}}}

\newcommand{\He}{{\ensuremath{\operatorname{Hess}\,}}}
\newcommand{\Tr}{{\ensuremath{\operatorname{Trace}}}}
\newcommand{\Cp}{{\ensuremath{C_{\rm pol}}}}
\newcommand{\Cpp}{{\ensuremath{C^1_{\rm pol}}}}

\newcommand{\V}{{\ensuremath{\mathcal{V}}}}

\newcommand{\Y}{{\ensuremath{\mathcal{Y}}}}
\newcommand{\cZ}{{\ensuremath{\mathcal{Z}}}}

\DeclareFontEncoding{LS1}{}{}
\DeclareFontSubstitution{LS1}{stix}{m}{n}
\DeclareMathAlphabet{\mathscr}{LS1}{stixscr}{m}{n}
\begin{document}
\title{Deep learning based  numerical  
approximation\\ algorithms for stochastic
partial differential equations 
}

\author{Christian Beck$^1$, Sebastian Becker$^1$, 
	Patrick Cheridito$^1$, \\
	Arnulf Jentzen$^{2,3}$ and Ariel Neufeld$^4$
	\bigskip
	\\
	\small{$^1$ Department of Mathematics, ETH Zurich}\\
	\small{$^2$ School of Data Science and School of Artificial Intelligence,}\\
	\small{The Chinese University of Hong Kong, Shenzhen (CUHK-Shenhen)}\\
	\small{$^3$ 
		Applied Mathematics: Institute for Analysis and Numerics, 
		University of Münster
	}\\
	\small{$^4$ Division of Mathematical Sciences, School of Physical}
	\\
	\small{and Mathematical Sciences, Nanyang Technological University}
}

\date{} 

\maketitle
\begin{abstract}
In this article, we introduce a deep learning based approximation algorithm for SPDEs.
Our approach employs neural networks to approximate the solutions of SPDEs along given 
realizations of the driving noise process. 
If applied to a set of simulated noise trajectories, it yields empirical distributions of 
SPDE solutions, from which functionals like the mean and variance can be estimated.
We test the performance of the method on stochastic heat equations with additive and multiplicative noise as well as stochastic Black–Scholes equations with multiplicative noise 
and Zakai equations from nonlinear filtering theory. In all cases, the proposed algorithm yields 
accurate results with short runtimes in up to 100 space dimensions.
\end{abstract}
\vspace{0.2cm}

\section{Introduction}
Stochastic partial differential equations (SPDEs) are key ingredients in a number 
of models in the natural sciences, engineering and finance. For example, they
are used to describe random surfaces in surface growth models \cite{HairerKPZ, BlomkerRomito2015},
temporal dynamics in Euclidean quantum field theories \cite{MourratWeber2015}, velocity fields in turbulent flows 
\cite{Birnir2013a, Birnir2013b}, prices of interest-rate derivatives \cite{FilipovicTappeTeichmann2010, HarmsStefanovitsTeichmannWutrich2015} or the spread of a biological or chemical contaminant in a river, a water basin or the groundwater system \cite{KallianpurXiong1994, KouritzinLong2002}.
SPDEs also play an important role in nonlinear filtering. For instance, the 
Kushner \cite{kushner1964differential} and Zakai equations \cite{zakai1969optimal}
describe the conditional density of a signal inferred from noisy observations; see, e.g.,
\cite{BudmanHolcombMorari91,ChenSun91,Rutzler87,SeinfeldGavalasHwang71,SolimanRay79b,WindesCinarRay89} for applications in chemical engineering, 
\cite{ duc2015ensemble,buehner2017ensemble,cassola2012wind,che2016wind,falissard2013genuinely,pelosi2017adaptive} for applications in weather forecasting and
\cite{BrigoHanzon98,CeciColaneri17,CoculescuGemanJeanblanc08,DuffieLando01,FreyRunggaldier10,FreySchmidt12} for applications in financial engineering.

Most SPDEs appearing in applications cannot be solved explicitly, and some of them 
are high-dimensional. Typical numerical approximation schemes discretize 
the time interval $[0,T]$ and the state space; see, e.g., \cite{debussche2011weak,debussche2009weak,du2002numerical,gyongy1999lattice,hausenblas2002numerical,hausenblas2003approximation,hausenblas2008finite,shardlow1999numerical,walsh2005finite,yan2005galerkin}
for temporal discretizations based on the linear implicit Euler method,
\cite{hutzenthaler2016strong,jentzen2009overcoming,lord2010stochastic,mukam2016note,wang2015note}
 for temporal discretizations based on exponential Euler-type  methods, 
 \cite{hausenblas2002numerical,hausenblas2003approximation,shardlow1999numerical,walsh2005finite}
for temporal discretizations based on linear implicit Crank--Nicolson-type methods,
\cite{bayer2016splitting,bensoussan1990approximation,bessaih2014splitting,florchinger1991time,gyongy2003splitting,legland1992splitting}
for temporal discretizations based on splitting-up approximation methods,
\cite{jentzen2011efficient,jentzen2015milstein,kruse2014consistency,reisinger2018stability, wang2014higher}
for temporal discretizations based on higher order approximation methods,
\cite{brzezniak2013finite,geissert2009rate,katsoulakis2011noise,kovacs2010finite,lord2010modified,walsh2005finite,yan2005galerkin}
for spatial discretizations based on finite elements methods,
\cite{gyongy2006numerical,JovanovicSuli,millet2005implicit,pettersson2005numerical,roth2006combination,shardlow1999numerical,walsh2006numerical}
for spatial discretizations based on finite differences methods and
\cite{flandoli2008introduction,grecksch1996time,hausenblas2003approximation,kloeden2001linear,lord2007postprocessing,muller2007implicit,mueller2008optimalOU}
for spatial discretizations based on spectral Galerkin methods.
Recently, deep neural networks have been used by \cite{zhang2020learning} to
approximate solutions to one-dimensional SPDEs. 
In \cite{teng2022solving}, a deep learning scheme for solving forward-backward doubly stochastic differential equations (FBDSDEs) has been derived, which by the stochastic version of the nonlinear Feynman--Kac formula \cite{pardoux1994backward}, is equivalent to solving certain semilinear backward SPDEs. In \cite{yao2021deep},  a deep learning scheme for a class of high-dimensional backward SPDEs with Neumann boundary condition has been developed by exploiting the extended stochastic nonlinear Feynman--Kac formula \cite{aman2013obstacle, Boufoussi2007}, which provides a stochastic representation for such SPDEs in terms of generalized FBDSDEs. \cite{teng2022solving,yao2021deep} both have numerical examples in
up to 100 space dimensions. Furthermore, a splitting scheme has been developed in \cite{bao2021solving} to decompose FBDSDEs into a BSDE and SDE, which, together with a Gauss--Hermite quadrature formula and Monte Carlo approximation, makes it possible to numerically solve one-dimensional FBDSDEs.
For more information on numerical methods for SPDEs we refer to 
the overview articles \cite{Gyongy02,jentzen2009numerical} 
and the monographs \cite{jentzen2011taylor,kruse2014strong}.
%

In this article we introduce and study a deep learning approximation method for 
high-dimensional SPDEs. It is inspired by the recent approximation methods for
nonlinear PDEs \cite{DeepSplitting,MR4293960,han2018solving} and 
uses deep neural networks to learn the solution of an SPDE along a given 
realization of the driving noise process. When applied to multiple simulated noise trajectories, it produces 
empirical distributions of SPDE solutions, enabling estimation of functionals like the mean and variance.
We test  the performance of the method in different numerical examples.
%
%
In each of them it produces accurate results with short runtimes in up to 100 space dimensions. 

The rest of this paper is organized as follows. In Section~\ref{sec:derivation} we derive 
the proposed approximation algorithm. In Section~\ref{sec:examples} we test its performance 
on stochastic heat equations with additive noise and multiplicative noise,
stochastic Black--Scholes equations with multiplicative noise and Zakai equations.

\section{Derivation of the proposed approximation algorithm
	\label{sec:derivation}
}
Let $ T \in (0,\infty) $ and $ d, \delta \in \N $. Consider a probability space
$(\Omega, \mathcal{F}, \P)$ equipped with a filtration 
$( \mathcal{F}_t )_{ t \in [0,T] }$ satisfying the usual conditions\footnote{
i.e. $( \mathcal{F}_t )_{ t \in [0,T] }$ is right-continuous and ${\mathcal{F}_0}$ 
contains all subsets of $\P$-null sets in $\mathcal{F}$} and let 
\[
Z=(Z_t(x,\omega))_{(t,x,\omega)\in [0,T]\times \R^d \times \Omega} 
\colon [0,T] \times \R^d \times \Omega \to \R^{\delta}
\]
be a sufficiently regular random field such that for every
$
x \in \R^d,
$ 
$
(Z_t(x))_{t\in [0,T]} \colon [0,T] \times \Omega \to \R^{\delta}
$ 
is an 
$
(\mathcal F_t)_{t \in [0,T]}$-It\^o process. Consider sufficiently regular functions
$
\varphi\colon \R^d \to \R,
$
$ 
f \colon \R^d \times \R \times \R^d \to \R,  
$
$
b \colon \R^d \times \R \times \R^d \to \R^{\delta},
$
$ 
\mu \colon \R^d \to \R^d,
$ 
$
\sigma \colon  \R^d \to \R^{ d \times d },
$
and let 
$
X=(X_t(x,\omega))_{(t,x,\omega)\in [0,T]\times \R^d \times \Omega}\colon [0,T]\times\R^d\times\Omega\to\R
$ 
be a random field such that for all
$
t\in [0,T]
$ and
$
x\in\R^d,
$ 
$
X_t(x)\colon\Omega\to\R
$ 
is 
$\mathcal F_t$/$\B(\R)$-measurable, for every $\omega\in\Omega$,
$
(X_t(x,\omega))_{(t,x)\in [0,T]\times\R^d}$ is in $C^{0,2}([0,T]\times\R^d,\R)
$ 
with at most polynomially growing partial derivatives of order 0, 1 and 2
with respect to the $x$-variables,
and for all $t\in [0,T]$ and $x\in\R^d$, one has
\begin{equation}\label{eq:SPDE}
\begin{split}
X_{ t }( x )
&  =
\varphi( x ) 
+
\int_{ 0 }^{ t }
f\big( 
x, X_s(x), \nabla X_s ( x ) 
\big)
\, ds
+
\int_{ 0 }^{ t } \big\langle b\big( x, X_s( x ), \nabla X_s ( x ) \big), dZ_s(x) \big\rangle_{\R^\delta}
\\
& \quad
+
\int_{ 0 }^{ t }
\Big(
\tfrac{ 1 }{ 2 } \Tr \brak{
\sigma( x ) \sigma(x )^T \He X_s ( x )} +
\big\langle \mu( x ), \nabla X_s ( x ) \big\rangle_{ \R^d } \Big) ds \; \; \; \mbox{$\P$-a.s.},
\end{split}
\end{equation}
where $\nabla$ and Hess denote the gradient and Hessian with respect to the $x$-variables,
respectively. Our goal is to compute an approximation of the solution  
$X\colon [0,T]\times\R^d \times \Omega \to \R$ of the SPDE \eqref{eq:SPDE}.
\subsection{Temporal discretization}
\label{subsec:temp-discret}
In this subsection we discretize the SPDE 
\eqref{eq:SPDE} in time by employing the splitting-up method 
(see, e.g., \cite{Bensoussan_SplittingUpMethodForSPDEs1992,
	BensoussanGlowinskiRascanu_ApproximationBySplittingUp1992,
	GrekschLisei_ApproximationOfStochasticNonlinearEquationsBySplittingMethod2013,
	gyongy2003rate,
	gyongy2003splitting,
	legland1992splitting}) 
to obtain a semi-discrete approximation problem. More precisely, we choose an $N\in\N$
and let $t_0, t_1, \ldots, t_N$ be real numbers such that
\[
\label{eq:time-step-discrete}
0 = t_0  < t_1 < \ldots < t_N = T.
\]
It follows from \eqref{eq:SPDE} that for all
$n\in\{0,1,\ldots,N-1\}$, $t \in (t_n,t_{n+1}]$ and $x\in\R^d$, we have
\[
\begin{split}
X_t(x) 
& 
= X_{t_n}(x) 
+
\int_{ t_n }^{ t }
f\big( 
x, X_s(x), \nabla X_s ( x ) 
\big)
\, ds
+
\int_{ t_n }^{ t }
\big\langle b\big( x, X_s( x ), \nabla X_s ( x ) \big), dZ_s(x) \big\rangle_{\R^\delta}
\\
& \quad
+ \int_{t_n}^t 
\Big(
\tfrac{ 1 }{ 2 }
\Tr \big( 
\sigma( x ) \sigma( x )^T \He X_s( x )
\big)
+
\big\langle \mu( x ), \nabla X_s ( x ) \big\rangle_{ \R^d }
\Big) ds \quad \mbox{$\P$-a.s.,}
\end{split}
\]
and therefore,
\[
\begin{split}
X_t (x) 
& \approx 
X_{t_n}(x) +
\int_{ t_n }^{ t_{ n + 1 } }
f\big( 
x, X_s(x), \nabla X_{s}( x ) 
\big)
\, ds
+
\int_{ t_n }^{ t_{ n + 1 } }
\big\langle b\big( x, X_s( x ), \nabla X_s ( x ) \big), dZ_s(x) \big \rangle_{\R^\delta}
\\
& \quad 
+ \int_{t_{n}}^t 
\Big(
\tfrac{ 1 }{ 2 } \Tr \big( 
\sigma( x ) \sigma( x )^T \He X_s ( x )
\big) 
+ 
\big\langle \mu( x ), \nabla X_s( x ) \big\rangle_{ \R^d }\Big) ds,
\end{split}
\]
which, in turn, suggests that 
\be \label{eq:Xapproximation}
\begin{split} 
X_{t}(x)  
&
\approx
{\cal H}_n \big( x, X_{t_n}(x), \nabla X_{t_n}(x) , Z_{t_{n+1}}(x) - Z_{t_n}(x) \big)
\\
& \quad
+ 
\int_{t_{n}}^t 
\Big(
\tfrac{ 1 }{ 2 }
\Tr \big( 
\sigma( x ) \sigma( x )^T \He X_s ( x )
\big) 
+ 
\big\langle \mu( x ), \nabla X_s ( x ) \big\rangle_{ \R^d } \Big) ds 
\end{split}
\ee
if ${\cal H}_n \colon \R^d \times \R \times \R^d \times \R^{\delta} \to \R$ are suitable
functions such that 
\beas
&&{\cal H}_n \big( x, X_{t_n}(x), \nabla X_{t_n}(x) , Z_{t_{n+1}}(x)-Z_{t_n}(x) \big)\\
&\approx& X_{t_n}(x) +
\int_{ t_n }^{ t_{ n + 1 } }
f\big( 
x, X_s(x), \nabla X_{s} ( x ) 
\big) ds +
\int_{ t_n }^{ t_{ n + 1 } }
\big\langle b\big( x, X_s( x ), \nabla X_s( x ) \big), dZ_s(x) \big\rangle_{\R^\delta}
\eeas
for all $n \in \{0, \dots, N-1\}$. For instance, in
the examples of Subsections \ref{subsec:stoch_heat} and \ref{subsec:Zakai} below we use the
first order approximation
\be \label{1storder}
\begin{aligned} 
&
{\cal H}_n \brak{x, X_{t_n}(x), \nabla X_{t_n}(x)\big) , \big(Z_{t_{n+1}}(x)-Z_{t_n}(x)}
= X_{t_n}(x) \, + \\
& \quad f
\bigl(
x,X_{t_n}(x), \nabla X_{t_n}(x)
\bigr)\,(t_{n+1}-t_n) + \big\langle
b\big(x, X_{t_n}(x), \nabla X_{t_n}(x)\big) , \big(Z_{t_{n+1}}(x)-Z_{t_n}(x)\big)
\big\rangle_{\R^{\delta}}.
\end{aligned}
\ee
But for the examples of Subsections \ref{subsec:const-coeff}--\ref{subsec:BS}, we obtained 
better numerical results for Milstein-type second order approximations; see 
Section \ref{sec:examples} below for more details.

Now, let
$
U \colon (0,T]\times\R^d\times\Omega\to\R
$ 
be a random field such that for all 
$
\omega \in \Omega
$ and
$
n \in \{0,1,\ldots,N-1\},
$
$
(U_t(x,\omega))_{(t,x)\in (t_n, t_{n+1}]\times\R^d}$ is in
$
C^{1,2}( (t_n, t_{n+1}]\times\R^d,\R)
$
with at most polynomially growing partial derivatives of order 0, 1 and 2
with respect to the $x$-variables, for all $\omega\in\Omega$ and $x\in\R^d$, one has
$
\int_0^T \big( \|\He U_s(x,\omega)\|_{\R^{d\times d}} + \|\nabla U_s(x,\omega)\|_{\R^d} \big) ds < \infty
$,
and for all $n \in \{0,1,\ldots,N-1\}$, $t \in (t_n,t_{n+1}]$ and $x \in \R^d$, 
\begin{equation}\label{eq:mildFormulationUPDE}
\begin{split}
U_t(x) 
& =
{\cal H}_n \big( x, X_{t_n}(x), \nabla X_{t_n}(x) , Z_{t_{n+1}}(x)-Z_{t_n}(x) \big)\\
& \quad
+  
\int_{t_n}^t \Big(
\tfrac{ 1 }{ 2 }
\Tr \big( 
\sigma( x ) \sigma( x )^T \He U_s( x )
\big) 
+ 
\big\langle \mu( x ), \nabla U_s( x ) \big\rangle_{ \R^d }
\Big) ds.
\end{split}
\end{equation}
\eqref{eq:Xapproximation} and \eqref{eq:mildFormulationUPDE} suggest that
$U_{t_{n}}(x) \approx X_{t_{n}}(x)$ for all $n\in\{1,2,\ldots,N\}$ and $x\in\R^d$.
Next, let $V \colon [0,T]\times\R^d\times\Omega\to\R$  
be a random field such that for all $x \in \R^d$, $V_0(x) = \varphi(x),$ 
for all $\omega\in\Omega$ and $n\in\{0,1,\ldots,N-1\}$,
$(V_t(x,\omega))_{(t,x)\in (t_n,t_{n+1}]\times\R^d}$
belongs to $C^{1,2}((t_n,t_{n+1}]\times\R^d,\R)$ with at most
polynomially growing partial derivatives of order 0, 1 and 2 with respect to the $x$-variables, for all
$
\omega \in \Omega
$ and
$
x \in \R^d,
$ one has
$
\int_0^T \big(\| \He V_s(x,\omega) \|_{\R^{d\times d}}
 + 
\|\nabla V_s(x,\omega)\|_{\R^d} \big) ds 
< \infty
$
and for all 
$
n \in \{0,1,\ldots,N-1\}
$,
$
t \in (t_n, t_{n+1}]
$ and
$
x \in \R^d,
$
\begin{equation}\label{eq:mildFormulationVPDE}
\begin{split}
V_t(x) 
& =
{\cal H}_n \big( x, V_{t_n}(x), \nabla V_{t_n}(x) , Z_{t_{n+1}}(x)-Z_{t_n}(x) \big)\\
& \quad + 
\int_{ t_n }^{ t }
\Big(
\tfrac{ 1 }{ 2 }
\Tr \big( 
\sigma( x ) \sigma( x )^T \He V_s( x )
\big)
+
\big\langle \mu( x ), \nabla V_s ( x ) \big\rangle_{ \R^d }
\Big)
\, ds
\end{split}
\end{equation}
(see, e.g., Deck \& Kruse \cite{DeckKruse_ParametrixMethod2002},
Hairer et al.\ \cite[Section 4.4]{HairerHutzenthalerJentzen_LossOfRegularity2015}, 
Krylov~\cite[Chapter 8]{Krylov_LecturesHoelder1996},
or
Krylov~\cite[Theorem 4.32]{Krylov_KolmogorovsEquations1998}
for existence, uniqueness and regularity results for 
\eqref{eq:mildFormulationUPDE} and \eqref{eq:mildFormulationVPDE}).
$V$ is a splitting-up type approximation of the 
random field $X$ (see, e.g., \cite{Bensoussan_SplittingUpMethodForSPDEs1992,
	BensoussanGlowinskiRascanu_ApproximationBySplittingUp1992,
	GrekschLisei_ApproximationOfStochasticNonlinearEquationsBySplittingMethod2013,
	gyongy2003rate, gyongy2003splitting, legland1992splitting}), and by
\eqref{eq:Xapproximation}--\eqref{eq:mildFormulationVPDE}, one has
\begin{equation} \label{eq:VapproxX}
	V_{t_{n}}(x) \approx U_{t_{n}}(x) \approx X_{t_{n}}(x) \quad \mbox{
		for all $n\in\{1,2,\ldots,N\}$ and $x\in\R^d$.}
\end{equation}

\subsection{Pathwise approximation}

Now, let us denote by $\Cp(\R^m, \R^n)$ the set of all continuous functions 
$\psi \colon \R^m \to \R^n$ that grow at most polynomially and by $\Cpp(\R^m, \R)$ 
the set of all differentiable functions $\psi \colon \R^m \to \R$ such that $\nabla \psi \in \Cp(\R^m, \R^m)$.
We assume that $\varphi \in \Cp(\R^m;\R)$,
\begin{equation} \label{Hn}
x \mapsto \mathcal{H}_n \bigl(x,v(x),\nabla v(x),w(x)\bigr)\in \Cp(\mathbb R^d;\mathbb R)
\end{equation}
for all $n\in\{0,1,\ldots,N-1\}$, $v \in \Cpp(\mathbb R^d;\mathbb R)$,
$w \in \Cp(\mathbb R^d;\mathbb R^\delta)$, and for every $\psi \in \Cp(\mathbb R^d;\mathbb R)$, 
there exists a unique function $u^\psi \in C^{1,2}((0,T]\times \R^d,\R)$ with at most polynomially 
growing partial derivatives of order 0, 1 and 2 with respect to the x-variables that depends 
continuously on $\psi$ such that
\begin{equation} \label{upsi}
u^\psi(t,x)= \psi(x)+\int_{0}^{ t }
\Big( \tfrac{ 1 }{ 2 }
\Tr \big( \sigma( x ) \sigma( x )^T \He u^\psi(s, x )
\big)
+
\big\langle \mu( x ), \nabla u^\psi(s, x ) \big\rangle_{ \R^d }
\Big)
\, ds
\end{equation}
for all $t \in [0,T]$ and $x \in \R^d$.

\begin{remark}
Clearly, \eqref{Hn} is satisfied if $\mathcal{H}_n$ is e.g. of the form 
\eqref{1storder} with $f, b \in \Cp(\R^d\times \R \times \R^d,\R)$, and
it is well known that \eqref{upsi} holds under standard assumptions 
such as smooth coefficients $\mu,\sigma$ with bounded derivatives and uniformly elliptic 
$\sigma\sigma^T$.
\end{remark}

For $\psi \in \Cpp$, we denote $P_t \psi(x) = u^\psi(t,x)$, and for 
$z \in \cZ := \Cp([0,T] \times \R^d ; \R^\delta)$, we define
\be \label{defV}
\cV_0^{(z)}(x) := \varphi(x) \quad \mbox{and} \quad
\cV_t^{(z)}(x) := P_{t-t_n} \mathcal{H}_n\bigl(x, \cV_{t_n}^{(z)}(x),\nabla 
\cV_{t_n}^{(z)}(x),z_{t_{n+1}}(x)-z_{t_n}(x)\bigr)
\ee
for $n\in\{0,1,\ldots,N-1\}$, $t\in (t_n,t_{n+1}]$ and $x\in \R^d$.	Then, for all $z \in \cZ$, 
$(\cV^{(z)}_t(x))_{(t,x) \in [0,T] \times \R^d}$ is uniquely given by $\cV^{(z)}_0(x) = \varphi(x)$ and
\begin{equation} \label{cV}
		 	\begin{split}
		 \mathcal{V}_t^{(z)}(x)
		 &=
		 \mathcal{H}_n\!\left(
		 x,
		 \mathcal{V}_{t_n}^{(z)}(x),
		 \nabla \mathcal{V}_{t_n}^{(z)}(x),
		 z_{t_{n+1}}(x)-z_{t_n}(x)
		 \right)\\
		 & \quad + 
		 \int_{ t_n }^{ t }
		 \Big(
		 \tfrac{ 1 }{ 2 }
		 \Tr \big( 
		 \sigma( x ) \sigma( x )^T \He \mathcal{V}^{(z)}_s( x )
		 \big)
		 +
		 \big\langle \mu( x ), \nabla \mathcal{V}^{(z)}_s ( x ) \big\rangle_{ \R^d }
		 \Big)
		 \, ds,
		\end{split}
		 \end {equation}
$n \in\{0,1,\ldots,N-1\}$, $t\in (t_n,t_{n+1}]$, $x\in \R^d$, and it follows 
from the definition \eqref{defV} that 
\be \label{meas}
\cV^{(z)}_t(x) \mbox{ is measurable in $(t,x,z) \in [0,T] \times \R^d \times \cZ$.}
\ee

%
\subsection{Approximate Feynman--Kac representation}
\label{subsec:approx_Feynman-Kac}

We now derive a Feynman--Kac representation for the pathwise splitting
solution \eqref{cV}. Let us assume that the driving noise process satisfies $\P(Z \in \cZ) = 1$.
By \eqref{meas}, $\mathcal{V}_t^{(Z(\omega))}(x)$ is measurable in 
$(t,x,\omega) \in [0,T] \times \R^d \times \Omega$, and for all $\omega \in \Omega$
such that $Z(\omega)\in\mathcal{Z}$, the function $\mathcal{V}_t^{(Z(\omega))}(x)$ solves equation~\eqref{eq:mildFormulationVPDE}. 
Consequently, if $V: [0,T]\times \R^d \times \Omega \to \R$ is a random field satisfying \eqref{eq:mildFormulationVPDE}, we obtain from $\omega$-wise uniqueness that

\be \label{V=cV}
V_t(x,\omega) = \mathcal{V}_t^{(Z(\omega))}(x) \quad \mbox{$\mathbb{P}$-a.s.}
\ee
%
Combining \eqref{V=cV} with the splitting approximation relation \eqref{eq:VapproxX}, we obtain
\[ \label{eq:Vz-approxX}
\mathcal V_{t_n}^{(Z(\omega))}(x)
= V_{t_n}(x,\omega)
\approx
X_{t_n}(x,\omega) 
\mbox{ for all } n\in \{0,1,\dots,N\}, x \in\mathbb R^d \mbox{ and $\P$-almost all } \omega \in \Omega. 
\]
Hence 
\be \label{condV}
\mathbb E[V_{t_n}(x)\mid Z]
=
\mathcal{V}_{t_n}^{(Z)}(x),
\ee
which suggests that 
\[ \label{eq:VapproxXWithZPluggedIn-b}
\V^{(z)}_{t_n}(x) \approx \E\big[ X_{t_n}(x)\!~|\!~Z=z\big] \quad
\mbox{for all $n \in \{0,1,\dots,N\}$, $x\in \R^d$ and $\mathbb{P}\circ Z^{-1}$-almost all $z\in \mathcal{Z}$.}
\]
%
In the following we introduce an auxiliary stochastic process allowing us to
apply a Feynman--Kac representation to \eqref{cV}. 
Assume $Y \colon [0, T] \times \Omega\to\R^d$ is a strong solution of
\begin{equation}\label{eq:SDEY}
Y_{t} = \xi + \int_0^t \mu(Y_s)\,ds 
+ \int_0^t \sigma(Y_s)\,dB_s,
\end{equation}
where $\xi\colon\Omega\to\R^d$ is an 
$\mathcal F_0$/$\B(\R^d)$-measurable initial 
condition satisfying $\E[\|\xi\|_{\R^d}^p] < \infty$ for all $p\in (0,\infty)$ and 
$B \colon [0,T]\times\Omega\to\R^d$
a standard $( \mathcal{F}_t )_{ t \in [0,T] }$-Brownian motion
such that $(\xi, B)$ is independent of $Z$. It is well known that if the
functions $\mu \colon \R^d \to \R^d$ and $\sigma \colon  \R^d \to \R^{ d \times d }$
grow at most linearly, one has
\begin{equation}
\label{eq:Y-integ-cond}
\sup_{t \in [0,T]} \E\big[\|Y_t\|_{\R^d}^p\big] < \infty
\quad \mbox{for all $p \in (0,\infty).$}
\end{equation}
From It\^{o}'s formula and the dynamics \eqref{eq:SDEY}, we obtain
\begin{equation} \label{eq:Ito}
\begin{split}
\V^{(z)}_{T-t}(Y_t) 
& = 
\V^{(z)}_{T-r}(Y_r) + 
\int_{r}^{t} 
\tfrac{\partial}{\partial s} 
\V^{(z)}_{T-s} 
(Y_s)\,ds + \int_{r}^{t} 
\big\langle \nabla \V^{(z)}_{T-s}(Y_s), \sigma(Y_s)\,dB_s  \big\rangle_{\R^d}\\
& \quad +
\int_{r}^{t} 
\tfrac{ 1 }{ 2 }
\operatorname{Trace}\!\big( 
\sigma(Y_s ) \sigma(Y_s )^T
\He \V^{(z)}_{T-s}( Y_s )
\big)\,ds \\
&  \quad + 
\int_{r}^{t}
\big\langle \mu( Y_s ), 
\nabla \V^{(z)}_{T-s} ( Y_s ) 
\big\rangle_{ \R^d }
\, ds \quad \mbox{$\P$-a.s.}
\end{split}
\end{equation}
for all $z \in \mathcal{Z}$, 
$n \in \{0, \dots, N-1\}$ and $r,t \in [T-t_{n+1}, T- t_n)$ such that $r < t$.
Moreover, it can be seen from \eqref{cV} that 
for all $n\in\{0,1,\ldots,N-1\}$ and $x\in\R^d$,
\[
\tfrac{\partial}{\partial t}  \V^{(z)}_{t}(x) 
= \big\langle \mu( x ), \nabla \V^{(z)}_{t}( x ) \big\rangle_{ \R^d } 
+ \tfrac{ 1 }{ 2 }
\Tr \big( 
\sigma( x ) \sigma( x )^T
\He \V^{(z)}_{t} ( x ), \quad t\in (t_{n},t_{n+1}),
\]
and therefore,
\[ \label{eq:backwardsEquationPoppingUpInIto}
\tfrac{\partial}{\partial t} \V^{(z)}_{T-t}(x) 
+ \big\langle \mu( x ), \nabla \V^{(z)}_{T-t} ( x ) \big\rangle_{ \R^d } 
+ \tfrac{ 1 }{ 2 }
\Tr \big( 
\sigma( x ) \sigma( x )^T
\He \V^{(z)}_{T-t} ( x )
\big) = 0, \quad t\in (T-t_{n+1},T-t_{n}).
\]
Combining this with \eqref{eq:Ito} gives 
\begin{equation}\label{eq:ItoOnOpenIntervalAfterCancellation}
\V^{(z)}_{T-t}(Y_t)  = \V^{(z)}_{T-r}(Y_r)
+ \int_{r}^t 
\big\langle \nabla \V^{(z)}_{T-s}(Y_s), 
\sigma(Y_s)\,dB_s  \big\rangle_{\R^d} \quad \mbox{$\P$-a.s.}
\end{equation}
for all $n \in \{0, \dots, N-1\}$ and $r,t \in [T-t_{n+1}, T- t_n)$ such that $r < t$.
If $\sigma \colon \R^d \to \R^{d\times d}$ and $\nabla \V^{(z)} \colon 
(t_n,t_{n+1}] \times \R^d \to \R^d$ are at most polynomially growing, one obtains 
from \eqref{eq:Y-integ-cond} that
\[
\int_{T-t_{n+1}}^t \E \Big[\big\| \sigma(Y_s)^T \nabla \V^{(z)}_{T-s}(Y_s)\big\|_{\R^d}^2\Big]\, 
ds <\infty, \quad t \in [T - t_{n+1}, T - t_n),
\] 
and as a consequence,
\[ \label{eq:conditionalExpectation1stTerm}
\E\bigg[\int_{T-t_{n+1}}^{t} 
\big\langle \nabla \V^{(z)}_{T-s}(Y_s), 
\sigma(Y_s)\,dB_{s} \big\rangle_{\R^d}\!~\Big|\!~\mathcal F_{T-t_{n+1}} \bigg] = 0,
\quad t \in [T - t_{n+1}, T - t_n).
\]
So, since $Y_{T-t_{n+1}}$ is 
$\mathcal{F}_{T- t_{n+1}}$/$\B(\R^d)$-measurable, one obtains from
\eqref{eq:ItoOnOpenIntervalAfterCancellation} that
\be \label{condY}
\E\big[ \V^{(z)}_{T-t}(Y_t)\!~\big|\!~ Y_{T-t_{n+1}}
\big]
= \E\big[
\V^{(z)}_{t_{n+1}}(Y_{T-t_{n+1}})\!~\big|\!~ Y_{T-t_{n+1}}
\big] = \V^{(z)}_{t_{n+1}}(Y_{T-t_{n+1}}) \quad \mbox{$\P$-a.s.}
\ee
for all $n\in\{0,1,\ldots,N-1\}$ and $t \in [T-t_{n+1}, T-t_{n})$.
Next, observe that since $(\V^{(z)}_t(x))_{(t,x)\in (t_n,t_{n+1}]\times\R^d}$ is in
$C^{1,2}((t_n,t_{n+1}]\times\R^d,\R)$ with at most polynomial growth
and $Y$ is $\P$-a.s. continuous, one has
\begin{equation} \label{25a}
\limsup_{t \uparrow T-t_n} 
\Big|\V^{(z)}_{T-t}(Y_t(\omega))-
\V^{(z)}_{T-t}(Y_{T-t_n}(\omega))\Big|=0
\end{equation}
for all $n\in\{0,1,\ldots,N-1\}$ and $\P$-almost all $\omega \in \Omega$.
Since \[
\int_0^T \big(
\|\He \V^{(z)}_s(x)\|_{\R^{d\times d}} 
+ 
\|\nabla \V^{(z)}_s(x)\|_{\R^d} \big)
\,ds 
< \infty \quad \mbox{for all $x \in \R^d$},
\]
we obtain from \eqref{cV} and \eqref{25a} that 
\be \label{eq:ptwConv}
\begin{aligned}
& 
\limsup_{t \uparrow T-t_n} 
\Big|\V^{(z)}_{T-t}(Y_t(\omega)) \\
& \quad - {\cal H}_n \big(Y_{T-t_n}(\omega), \V^{(z)}_{t_n}(Y_{T-t_n}(\omega)), 
\nabla \V^{(z)}_{t_n}(Y_{T-t_n}(\omega)) , z_{t_{n+1}}(Y_{T-t_n}(\omega))
-z _{t_n}(Y_{T-t_n}(\omega)) \big) \Big| = 0
\end{aligned}
\ee
for all $n\in\{0,1,\ldots,N-1\}$ and $\P$-almost all $\omega\in\Omega$.
In addition, it follows from \eqref{eq:Y-integ-cond} together with our assumptions 
that $\V^{(z)}_0(x)=\varphi(x)$ for all $x \in \R^d$ and
$(\V^{(z)}_t(x))_{(t,x)\in (t_n,t_{n+1}]\times\R^d}\in C^{1,2}((t_n,t_{n+1}]\times\R^d,\R)$
has at most polynomially growing partial derivatives of order 0, 1 and 2
with respect to the $x$-variables that 
\begin{equation}\label{NRX}
\sup_{t \in [0,T]} \E\big[\|Y_t\|^p_{\R^d}\big] 
+
\sup_{t \in [0,T]} \E\big[|\V^{(z)}_{T-t}(Y_t)|^p\big] 
+
\sup_{t \in [0,T]} \E\big[\|\nabla\V^{(z)}_{T-t}(Y_{t})\|^p_{\R^d}\big]  <\infty
\end{equation}
for all $p \in (0, \infty)$. So for sufficiently regular 
$ {\cal H}_n \colon \R^d \times \R \times \R^d \times \R^{\delta} \to \R$, we have
\begin{equation}\label{NRZ}
\begin{split}
&\sup_{t \in [0,T]} 
\E\Big[
\big|\V^{(z)}_{t}(Y_{T-t})\big|^p\Big] + \\
&
\max_{n \in \{0,1,\dots,N-1\}}\E\Big[\Big|
{\cal H}_n \big(Y_{T-t_n}, \V^{(z)}_{t_n}(Y_{T-t_n}), 
\nabla \V^{(z)}_{t_n} (Y_{T-t_n}) , z_{t_{n+1}}(Y_{T-t_n})
-z _{t_n}(Y_{T-t_n}) \big)  \Big|^p
\Big] <\infty
\end{split}
\end{equation}
for all $p \in (0, \infty)$, which together with \eqref{eq:ptwConv}, yields 
\begin{align*}
& \limsup_{t\uparrow T-t_n} \E\Big[
\Big|
\V^{(z)}_{T-t}(Y_t) - \\
& 
{\cal H}_n \big(Y_{T-t_n}, \V^{(z)}_{t_n}(Y_{T-t_n}), 
\nabla \V^{(z)}_{t_n}(Y_{T-t_n}) , z_{t_{n+1}}(Y_{T-t_n})
-z _{t_n}(Y_{T-t_n}) \big) 
\Big|
\Big] 
= 0
\end{align*}
for all $n\in\{0,1,\ldots,N-1\}$. Therefore, we obtain from \eqref{condY} that 
for all $n\in\{0,1,\ldots,N-1\}$, 
\begin{equation}\label{eq:FC}
\begin{aligned}
& \V^{(z)}_{t_{n+1}}(Y_{T-t_{n+1}}) = \\ &
\E\Bigl[ 
{\cal H}_n \big(Y_{T-t_n}, \V^{(z)}_{t_n}(Y_{T-t_n}), 
\nabla \V^{(z)}_{t_n}(Y_{T-t_n}) , z_{t_{n+1}}(Y_{T-t_n})
-z _{t_n}(Y_{T-t_n}) \big)
\!~\big|\!~ Y_{T-t_{n+1}}
\Bigr] 
\end{aligned}
\end{equation}
$\P$-a.s., which is the Feynman--Kac representation we have been aiming for. 

%
\subsection{Formulation as recursive minimization problems}
In this subsection we reformulate the conditional expectations in
\eqref{eq:FC} as recursive minimization problems. 
Note that by \eqref{NRZ}, we have for all $z \in \mathcal{Z}$
and ${\cal H}_n \colon \R^d \times \R \times \R^d \times \R^{\delta} \to \R$, 
\[
\begin{aligned}
	& \E\Big[\big|{\cal H}_n \big(Y_{T-t_n}, \V^{(z)}_{t_n}(Y_{T-t_n}), 
	\nabla \V^{(z)}_{t_n}(Y_{T-t_n}) , z_{t_{n+1}}(Y_{T-t_n}) -z _{t_n}
	(Y_{T-t_n}) \big) \big|^2\Big] < \infty, \quad n \in \{0, \dots, N-1\}.
\end{aligned}
\]
Therefore, we obtain from \eqref{eq:FC} by the $L^2$-minimality property of conditional 
expectations (see, e.g., Klenke \cite[Corollary 8.17]{Klenke_2014}) and the 
factorizaton lemma (see, e.g., Klenke \cite[Corollary 1.97]{Klenke_2014})
that for all $z \in \mathcal{Z}$ and $n \in \{1,2,\dots,N\}$,
\begin{equation}\label{eq:formulationAsMinimizationProblem-L2}
	\begin{split} 
		&\V^{(z)}_{t_{n}}
		= 
		\operatornamewithlimits{argmin}_{u \in L^2(\P\circ Y^{-1}_{T-t_n})}
		\E
		\Big[ 
		\big|u(Y_{T-t_{n}}) - \\ &
		{\cal H}_{n-1} \big(Y_{T-t_{n-1}}, \V^{(z)}_{t_{n-1}}(Y_{T-t_{n-1}}), 
		\nabla \V^{(z)}_{t_{n-1}}(Y_{T-t_{n-1}}) , z_{t_{n}}(Y_{T-t_{n-1}}) -z _{t_{n-1}}(Y_{T-t_{n-1}}) \big)
		\big|^2
		\Big],
	\end{split}
\end{equation}
where the equality is understood in the $\P\circ Y^{-1}_{T-t_n}$-a.s. sense and
the minimizer is unique as an element of $L^2(\P\circ Y^{-1}_{T-t_n})$. Moreover, since we know that
$\cV^{(z)}_{t_n} \colon \R^d \to \R$ is continuous, it is enough to minimize over continuous functions 
$u \in L^2(\P\circ Y^{-1}_{T-t_n})$. 
%
\subsection{Deep neural network approximations}
In this subsection we approximate the functions $\V^{(z)}_{t_n} \colon \R^d \to \R$, 
$n\in\{1,2,\ldots,N\}$, with neural networks. More precisely, we choose $\nu\in\N$ and 
consider functions
$\bV_{n}=(\bV_{n}(\theta,x))_{(\theta,x)\in \R^\nu \times \R^d}\colon \R^{\nu}\times\R^d\to\R$, $n\in\{0,1,\dots,N\}$, 
%
such that $\bV_0(\theta,x) = \varphi(x)$ for all $\theta \in \R^{\nu}$
and $x\in\R^d$. For $n \in \{1, \dots, N\}$, we try to find a parameter vector 
$\theta^{z,n} \in \R^{\nu}$ such that 
\be
\label{eq:VThetaApproxVz}
\bV_n(\theta^{z,n} ,x) \approx \V^{(z)}_{t_n}(x) \quad \mbox{and therefore,} \quad
\bV_{n}(\theta^{z,n} ,x) \approx 
\E\big[X_{t_n}(x)\!~|\!~Z=z\big] \quad \mbox{for all } x \in \R^d.
\ee
%
%
We specify the functions $\bV_n\colon\R^{\nu}\times\R^d\to\R$, 
$n\in\{1,2,\ldots,N\}$, as neural networks.
%
%
(see, e.g, \cite{Bengio09,LeCunBengioHinton15}). For $k\in\N$, let 
$ \mathcal{T}_{k} \colon \R^k \to \R^k $ 
be given by 
\begin{equation}
\label{eq:logistic}
\mathcal{T}_k( x ) 
=
\big(
\!
\tanh(x_1), 
\tanh(x_2), 
\dots
,
\tanh(x_k) 
\big), \quad x \in \R^k.
\end{equation}
For $ \theta = ( \theta_1, \theta_2, \dots, \theta_{ \nu } ) \in \R^{ \nu } $, 
$ v \in \N_0 = \{0\} \cup \N$,
$ k, l \in \N $
with 
$
v + lk + l \leq \nu,
$
let 
$ A^{ \theta, v }_{ k, l } \colon \R^k \to \R^l $ be of the form
\begin{equation}
A^{ \theta, v }_{ k, l }( x )
=
\left(
\begin{array}{cccc}
\theta_{ v + 1 }
&
\theta_{ v + 2 }
&
\dots
&
\theta_{ v + k }
\\
\theta_{ v + k + 1 }
&
\theta_{ v + k + 2 }
&
\dots
&
\theta_{ v + 2 k }
\\
\theta_{ v + 2 k + 1 }
&
\theta_{ v + 2 k + 2 }
&
\dots
&
\theta_{ v + 3 k }
\\
\vdots
&
\vdots
&
\vdots
&
\vdots
\\
\theta_{ v + ( l - 1 ) k + 1 }
&
\theta_{ v + ( l - 1 ) k + 2 }
&
\dots
&
\theta_{ v +  lk }
\end{array}
\right)
\left(
\begin{array}{c}
x_1
\\
x_2
\\
x_3
\\
\vdots 
\\
x_k
\end{array}
\right)
+
\left(
\begin{array}{c}
\theta_{ v +  lk + 1 }
\\
\theta_{ v + lk + 2 }
\\
\theta_{ v + lk + 3 }
\\
\vdots 
\\
\theta_{ v + lk + l }
\end{array}
\right), \quad x \in \R^k.
\end{equation}
If $s\in \N$ satisfies $(sd+ 1)(d+1) \leq \nu$, one can define
the functions $\bV_n \colon\R^{\nu} \times\R^d \to \R$, $n\in\{1,\ldots,N\}$, by 
\begin{align}\label{def:NN}
& \bV_{ n }(\theta,x)
= A^{ \theta, sd(d+1) }_{ d, 1 } 
\circ 
\mathcal{T}_d
\circ
A^{ \theta, (s-1)d(d+1) }_{ d, d } 
\circ
\ldots
\circ 
\mathcal{T}_d
\circ 
A^{ \theta, d(d+1) }_{ d, d } 
\circ 
\mathcal{T}_d
\circ 
A^{ \theta, 0 }_{ d, d }(x), \quad \mbox{$\theta \in \R^{\nu}$, $x \in \R^d$.}
\end{align}
For every $n \in \{1,2,\dots,N\}$ ,
$\bV_n\colon\R^{\nu}\times\R^d\to\R$ describes a fully-connected feedforward 
deep neural network with $s+2$ layers (one input layer with $d$ neurons, $s$ 
hidden layers with $d$ neurons each, and one output layer with one neuron) and the hyperbolic tangent 
as activation function.


\subsection{Stochastic gradient descent based minimization}
In the next step, we attempt to find 
suitable parameter vectors $\theta^{z,1},\theta^{z,2}$, 
$\ldots,\theta^{z,N}\in\R^{\nu}$  for \eqref{eq:VThetaApproxVz}
by recursive minimization. More precisely, for given 
$n\in\{1,2,\ldots,N\}$ and $\theta^{z,0},\theta^{z,1},\dots,\theta^{z,n-1} \in \R^\nu$,
we employ a stochastic gradient descent algorithm to obtain an approximate minimizer
$\theta^{z,n}\in\R^{\nu}$ of the loss function
\begin{align} \label{eq:toMinimize}
& 
\theta
\mapsto 
\E\Big[\big| \bV_{n}(\theta,Y_{T-t_{n}}) - \\ \nonumber
&
{\cal H}_{n-1} \big(Y_{T-t_{n-1}}, 
\bV_{n-1}(\theta^{z,n-1},Y_{T-t_{n-1}}) , 
\nabla_x\bV_{n-1}(\theta^{z,n-1},Y_{T-t_{n-1}}),
z_{t_{n}}(Y_{T-t_{n-1}})-z_{t_{n-1}}(Y_{T-t_{n-1}})
\big)
\big|^2\Big].
\end{align}
This can be implemented by fixing a large $M \in \N$ and 
introducing for every $m \in \{0, \dots, M-1\}$, 
a strong solution $Y^{m}\colon [0,T]\times\Omega\to\R^d$ of the SDE
\begin{equation}
\label{eq:SDE-Y-m}
Y^{m}_t = \xi^{m} 
+ \int_0^t \mu(Y^{m}_s)\,ds 
+ \int_0^t \sigma(Y^{m}_s)\,dB^{m}_s, 
\end{equation}
where $\xi^{m}\colon \Omega \to \R^d$, $m\in \{0, \dots, M-1\}$, are 
i.i.d.~$\mathcal F_0/\B(\R^d)$-measurable functions and 
$B^{m} \colon [0,T] \times \Omega \to \R^d$, $m\in \{0, \dots, M-1\}$, 
are i.i.d.~standard 
$( \mathcal{F}_t )_{ t \in [0,T] }$-Brownian motions. Now, for a sufficiently 
small $\gamma\in (0,\infty)$, let
$\vartheta^{z,n}=(\vartheta^{z,n}_m)_{m = 0}^M \colon 
\{0, 1, \dots, M\} \times\Omega\to\R^\nu$, 
$n\in\{1,\ldots,N\}$, be stochastic processes satisfying 
\begin{align} \label{eq:sgdWithOriginalY}
&\vartheta^{z,n}_{m+1} 
=  
\vartheta^{z,n}_m 
- 
2\gamma\cdot
\nabla_{\theta}\bV_n(\vartheta^{z,n}_m,Y^{z,m}_{T-t_n})
\cdot 
\Big[
\bV_n(\vartheta^{z,n}_m,Y^{z,m}_{T-t_n}) 
- \\ &
{\cal H}_{n-1} \big(Y^{z,m}_{T-t_{n-1}}, 
\bV_{n-1}(\vartheta^{z,n-1}_M,Y^{z,m}_{T-t_{n-1}}),
\nabla_x\bV_{n-1}(\vartheta^{z,n-1}_M,Y^{z,m}_{T-t_{n-1}}),
z_{t_{n}}(Y^{z,m}_{T-t_{n-1}})-z_{t_{n-1}}(Y^{z,m}_{T-t_{n-1}})\big)
\nonumber 
\end{align}
for all $n \in \{1, \dots, N\} \mbox{ and } m \in \{0, \dots, M-1\}$.
%
%
\subsection{Temporal discretization of the auxiliary stochastic process}
\label{subsec:temp_disc}

Equation \eqref{eq:sgdWithOriginalY} provides an implementable numerical algorithm in 
the special case where the solution processes $Y^{m}\colon[0,T] \times \Omega \to \R^d$, 
$m\in \{0, \dots, M-1\}$, of the SDEs \eqref{eq:SDE-Y-m} can exactly be simulated.
If this is not possible, one can employ a numerical approximation method for SDEs, 
such as the Euler--Maruyama or Milstein scheme (see, e.g., Kloeden \& Platen
\cite{KloedenPlaten1992}). For instance, note that one obtains from
\eqref{eq:SDE-Y-m} for all $m \in \{0, \dots, M-1\}$ and $r,t \in [0,T]$ such that $r < t$,
\[
Y^{m}_{t} 
= 
Y^{m}_r 
+ 
\int_r^t \mu(Y^{m}_s)\,ds
+
\int_r^t \sigma(Y^{m}_s)\,dB^{m}_s \quad \mbox{$\P$-a.s.},
\]
and therefore, for all $n \in \{0, 1, \dots, N-1\}$,
\[
Y^{m}_{T-t_{n}} 
= 
Y^{m}_{T-t_{n+1}} 
+ 
\int_{T-t_{n+1}}^{T-t_n} \mu(Y^{m}_s)\,ds
+
\int_{T-t_{n+1}}^{T-t_n} \sigma(Y^{m}_s)\,dB^{m}_s \quad 
\mbox{$\P$-a.s.,}
\]
or written differently,
\[
Y^{m}_{\tau_{n+1}} = Y^{m}_{\tau_{n}}  + \int_{\tau_n}^{\tau_{n+1}} \mu(Y^{m}_s)\,ds  + \int_{\tau_n}^{\tau_{n+1}} \sigma(Y^{m}_s)\,dB^{m}_s \quad \mbox{$\P$-a.s.}
\]
for $0 = \tau_0 < \tau_1 < \dots < \tau_N = T$ given by $\tau_n=T-t_{N-n}$. This 
suggests that for all $n \in \{0,1, \dots, N-1\}$ and $m \in \{0, \dots, M-1\}$,
\[
Y^{m}_{\tau_{n+1}} \approx Y^{m}_{\tau_{n}}  + \mu( Y^{m}_{\tau_{n}}) \,(\tau_{n+1}-\tau_n)
+ \sigma( Y^{m}_{\tau_{n}})\, (B^{m}_{\tau_{n+1}}-B^{m}_{\tau_n}),
\]
leading to the Euler--Maruyama approximations
$\mathcal{Y}^{m} \colon \{0,1,\dots,N\} \times \Omega \to \R^d$ given by
\[
\label{EM} \Y^{m}_{0} = \xi^{m} \quad \mbox{and} \quad
\Y^{m}_{n+1} 
= 
\Y^{m}_{n} 
+ 
\mu(\Y^{m}_{n})\,(\tau_{n+1}-\tau_n)
+ 
\sigma(\Y^{m}_{n})\,(B^{m}_{\tau_{n+1}} - B^{m}_{\tau_n}). 
\]
Then 
\[
\Y^{m}_n\approx Y^{m}_{\tau_n}= Y^{m}_{T-{t_{N-n}}} \quad \mbox{or equivalently,}
\quad Y^{m}_{T-{t_{n}}} \approx \Y^{m}_{N-n} \]
for all $n \in \{0,1, \dots, N\}$ and $m \in \{0, \dots, M-1\}$,
which we use to approximate the stochastic processes $\vartheta^{z,n}
=(\vartheta^{z,n}_m)_{m = 0}^M \colon \{0, 1, \dots, M\} \times\Omega\to\R^\nu$, 
$n\in\{1,\ldots,N\}$, from \eqref{eq:sgdWithOriginalY} with processes
$\Theta^{z,n}=(\Theta^{z,n}_m)_{m = 0}^M \colon \{0, 1, \dots, M\} \times\Omega\to\R^\nu$, 
$n \in \{1,\ldots,N\}$, satisfying
\begin{align*}
&\Theta^{z,n}_{m+1} 
= 
\Theta^{z,n}_m 
- 
2\gamma\cdot
\nabla_{\theta}\bV_n(\Theta^{z,n}_m,\Y^{m}_{N-n})
\cdot
\Big[
\bV_n(\Theta^{z,n}_m,\Y^{m}_{N-n}) - \\ &
{\cal H}_{n-1} \big(\Y^{m}_{N-n+1}, 
\bV_{n-1}(\Theta^{z,n-1}_M,\Y^{m}_{N-n+1}) ,
\nabla_x\bV_{n-1}(\Theta^{z,n-1}_M,\Y^{m}_{N-n+1}),
\big(z_{t_{n}}(\Y^{m}_{N-n+1})-z_{t_{n-1}}(\Y^{m}_{N-n+1}) \big)
\Big].
\end{align*}
Then, one has 
\[
\Theta^{z,n}_M \approx \vartheta^{z,n}_M \quad \mbox{for all $n \in \{1, \dots, N\}.$}
\]
In the following two subsections we merge the derivations above 
and provide precise formulations of the proposed approximation algorithm; 
first in a special case in Subsection~\ref{subsec:algo1}) and
then in the general case in Subsection~\ref{subsec:algo-Full-gen}.
%
%
%
\subsection{Description of the proposed approximation algorithm in a special case}
\label{subsec:algo1}
In this subsection we describe the proposed approximation algorithm in the special case where 
the SDE \eqref{eq:SDE-Y-m} is discretized with the 
standard Euler--Maruyama scheme (see, e.g., Kloeden \& Platen \cite{KloedenPlaten1992} or
Maruyama \cite{Maruyama1955}), the approximation mappings 
${\cal H}_n \colon \R^d \times \R \times \R^d \times \R^{\delta}
\to \R$, $n \in \{0, \dots, N-1\}$ are the first order approximations \eqref{1storder} and
the minimization algorithm is plain vanilla stochastic gradient descent with a constant 
learning rate $\gamma \in (0,\infty)$ and batch size 1.
A description of a more general  algorithm allowing 
to incorporate more general SDE approximation schemes, more general 
approximation mappings ${\cal H}_n \colon \R^d \times \R \times \R^d \times \R^{\delta}
\to \R$, $n \in \{0, \dots, N-1\}$ and more sophisticated machine learning methods 
such as batch normalization (cf., for instance, Ioffe \& Szegedy \cite{IoffeSzegedy2015}) 
and Adam updating (cf., for example, Kingma \& Ba \cite{KingmaBa2015}), is given in Subsection~\ref{subsec:algo-Full-gen} below. 

\begin{algo}[Special case]
\label{algo:1}
	Let 
	$T,\gamma\in (0,\infty)$, 
	$d,\delta,N,M, s \in\N$ and $\nu = (sd+1)(d+1)$. Consider 
	$t_0,t_1,\ldots,t_N\in [0,T]$ such that
	$0 = t_0 < t_1 < \ldots < t_N = T$, and let
	$0 < \tau_0 < \tau_1 < \dots < \tau_n = T$ be given by $\tau_n= T-t_{N-n}$,
	$n \in \{0,1,\dots,N\}$. Let
	$
	f\colon \R^d\times \R\times \R^{d} \to \R
	$, 
	$
	b\colon \R^d\times \R\times \R^{d} \to \R^{\delta}
	$,
	$
	\mu\colon\R^d\to\R^d
	$
	and
	$
	\sigma\colon\R^d\to\R^{d\times d}
	$ be continuous functions and
	$(\Omega,\F,\P)$ a probability space equipped with a filtration 
	$({\cal F}_t)_{t \in [0,T]}$ satisfying the usual conditions.
	Consider i.i.d.\ ${\cal F}_0/{\cal B}(\R^d)$-measurable random variables $\xi^m$
	together with i.i.d.\ standard $({\cal F}_t)_{t \in [0,T]}$-Brownian motions
	$B^{m}\colon [0,T]\times\Omega\to\R^d$, $m\in \{0, \dots, M-1\}$. Consider stochastic 
	processes $\Y^{m}\colon \{0,1,\ldots,N\}\times\Omega\to\R^d$
		given by 
	\[
	\Y^{m}_0 = \xi^{m} \quad \mbox{and} \quad
	\Y^{m}_{n+1} 
	= 
	\Y^{m}_n 
	+ 
	\mu(\Y^{m}_n)\,(\tau_{n+1}-\tau_{n})
	+ 
	\sigma(\Y^{m}_n)\,(B^{m}_{\tau_{n+1}}-B^{m}_{\tau_{n}}), 
	\]
	$n \in \{0, \dots, N-1\},$
	$m \in \{0, \dots, M-1\}$.
	Let $ \mathcal{T}_d \colon \R^d \to \R^d $ be given by
	\begin{equation}
	\label{eq:activation}
\mathcal{T}_d( x ) 
=
\big(
\!
\tanh(x_1), 
\tanh(x_2), 
\dots
,
\tanh(x_d) 
\big), \quad x \in \R^d.
	\end{equation}
For all $ \theta = ( \theta_1, \theta_2, \dots, \theta_{ \nu } ) \in \R^{ \nu }$,
	$ k, l \in \N $,
	$ v \in \N_0 = \{0\} \cup \N $
	with 
	$
	v + kl + l \leq \nu,
	$
	let 
	$ A^{ \theta, v }_{ k, l } \colon \R^k \to \R^l $ be given by
	\begin{equation}
	A^{ \theta, v }_{ k, l }( x ) = \bigg(\textstyle\sum\limits_{i=1}^ k x_i\, \theta_{v+i} + \theta_{v+kl+1}, \dots, \textstyle\sum\limits_{i=1}^ k x_i\, \theta_{v+(l-1)k+i} + \theta_{v+kl+l}  \bigg).
	\end{equation}
Let $\bV_n\colon\R^{\nu}\times\R^d\to\R$, $n\in\{0,1,\ldots,N\}$, be given by
	$\bV_0(\theta,x) = \varphi(x)$
	and  
	\[
		\bV_{ n }(\theta,x)
	=   A^{ \theta, sd(d+1) }_{ d, 1 } 
	\circ \mathcal{T}_d \circ  \dots
	\circ 
	\mathcal{T}_d
	\circ 
	A^{ \theta, d(d+1)}_{ d, d } 
	\circ 
	\mathcal{T}_d
	\circ 
	A^{ \theta, 0 }_{ d, d }(x) 
	\]
	for $n \in \{1, \dots, N\}$, $\theta \in \R^\nu$ and $x \in \R^d$.
	For a function $z \in \cZ$,
	let $\Theta^{z,n}\colon \{0, \dots, M\} \times\Omega\to\R^{\nu}$, $n \in \{1,2,\ldots,N\}$, 
	be  stochastic processes and 
	$\phi^{z,n,m} \colon \R^{\nu}\times\Omega\to\R$, 
	$\Phi^{z,n,m}\colon\R^{\nu}\times\Omega\to\R^{\nu}$, $n \in \{1, \dots, N\}$, 
	$m \in \{0, \dots, M-1\}$, mappings
	such that
	\begin{equation}\label{phi-z-n-m-Special}
	\begin{split}
	&\phi^{z,n,m}(\theta,\omega) 
	= \Big[ \bV_n\big(\theta,\Y^{m}_{N-n}(\omega)\big) 
	- \bV_{n-1}(\Theta^{z,n-1}_M(\omega),\Y^{m}_{N-n+1}(\omega)) 
	\, - (t_{n}-t_{n-1}) \\
	& \cdot f\big(\Y^{m}_{N-n+1}(\omega),\bV_{n-1}(\Theta^{z,n-1}_M(\omega),\Y^{m}_{N-n+1}(\omega)),
	\nabla_x\bV_{n-1}(\Theta^{z,n-1}_M(\omega),\Y^{m}_{N-n+1}(\omega))\big)\\
	& - b\big(\Y^{m}_{N-n+1}(\omega),\bV_{n-1}(\Theta^{z,n-1}_M(\omega),\Y^{m}_{N-n+1}(\omega)),
	\nabla_x\bV_{n-1}(\Theta^{z,n-1}_M(\omega),\Y^{m}_{N-n+1}(\omega))\big)\\
	& \cdot \big(z_{t_{n}}(\Y^{m}_{N-n+1}(\omega))-z_{t_{n-1}}(\Y^{m}_{N-n+1}(\omega))\big)
	\Big]^2, 
	\end{split}
	\end{equation}
	\begin{equation}
	\label{eq:plain-vanilla-SGD}
	\Phi^{z,n,m}(\theta,\omega) = \nabla_{\theta}\phi^{z,n,m}(\theta,\omega) 
	\quad \mbox{and} \quad
	\Theta^{z,n}_{m+1} = \Theta^{z,n}_m - \gamma\cdot\Phi^{z,n,m}(\Theta^{z,n}_m)
	\end{equation}
for all $n \in \{1, \dots, N\}$, $m \in \{0, \dots, M-1\}$, $\theta \in \R^\nu$ and $\omega \in \Omega$.
\end{algo}
Note that in the setting of Framework~\ref{algo:1} one has 
\begin{align}\label{eq:sgdWithApproximatedY-SPECIAL}
&\Theta^{z,n}_{m+1} 
= 
\Theta^{z,n}_m 
- 
2 \gamma\cdot
\nabla_{\theta}\bV_n(\Theta^{z,n}_m,\Y^{m}_{N-n})
\cdot
\Big[
\bV_n(\Theta^{z,n}_m,\Y^{m}_{N-n}) 
- \bV_{n-1}(\Theta^{z,n-1}_M,\Y^{m}_{N-n+1}) \nonumber\\ 
& \quad 
-f \big(\Y^{m}_{N-n+1},\bV_{n-1}(\Theta^{z,n-1}_M,\Y^{m}_{N-n+1}),
\nabla_x\bV_{n-1}(\Theta^{z,n-1}_M,\Y^{m}_{N-n+1})\big)\,(t_{n}-t_{n-1})\nonumber\\
&  \quad
-b\big(\Y^{m}_{N-n+1},\bV_{n-1}(\Theta^{z,n-1}_M,\Y^{m}_{N-n+1}),
\nabla_x\bV_{n-1}(\Theta^{z,n-1}_M,\Y^{m}_{N-n+1})\big)\\
&  \quad
\cdot \big(z_{t_{n}}(\Y^{m}_{N-n+1})-z_{t_{n-1}}(\Y^{m}_{N-n+1})\big)
\Big] \nonumber
\end{align}
for all $n \in \{1, \dots, N\}$ and $m \in \{0, \dots, M-1\}$. For sufficiently large 
$N,M \in \N$ and sufficiently small $\gamma \in (0,\infty)$, we think of 
$\mathbb V_n(\Theta^{z,n}_M,x) \colon \Omega \to \R$ 
as a suitable approximation 
\begin{equation}\label{V-Algo-result-Special}
\mathbb V_n(\Theta^{z,n}_M,x) \approx \E\big[X_{t_n}(x)\!~|\!~Z=z\big],
\end{equation}
where
$
X \colon [0,T]\times\R^d\times\Omega\to\R
$ 
is a solution of SPDE \eqref{eq:SPDE}.

\subsection{Description of the proposed approximation algorithm in the general case}
\label{subsec:algo-Full-gen}
In Framework~\ref{def:general_algorithm} below we present 
a general approximation algorithm which includes the algorithm of
Framework~\ref{algo:1} above as a special case. Compared to Framework~\ref{algo:1}, Framework~\ref{def:general_algorithm} also covers other SDE approximation schemes 
than Euler--Maruyama, higher order approximation 
mappings ${\cal H}_n \colon \R^d \times \R \times \R^d \times \R^{\delta} \to \R$,
$n \in \{0, \dots, N-1\}$, such as, e.g., Milstein-type second order approximations,
as well as more advanced machine learning methods, such as 
batch normalization (c.f. Ioffe \& Szegedy \cite{IoffeSzegedy2015})
and Adam updating (cf.\ Kingma \& Ba \cite{KingmaBa2015}).
%
%
\begin{algo}[General case]
	\label{def:general_algorithm}
	Let $T \in (0,\infty)$,  $M,N, d, \delta, \varrho, \nu, \varsigma \in \N$, 
	and $(J_m)_{m \in \N_0} \subseteq \N$. Consider $t_0,t_1,\ldots,t_N\in [0,T]$ 
	such that $0 = t_0 < t_1 < \ldots < t_N = T$ and let $0 = \tau_0 < \tau_1 < \ldots < \tau_N = T$
	be given by $\tau_n= T-t_{N-n}$, $n \in \{0,\dots,N\}$. Consider 
	functions
	$H\colon [0,T]^2\times\R^d\times\R^d\to\R^d$
	and
	$\mathcal{H}_n\colon \R^d \times \R \times \R^d \times \R^\delta \to \R$, $n\in \{0,\dots,N-1\}$.
	%
	Let 
	$ 
	(\Omega, \F, \P ) 
	$ 
	be a  probability space equipped with a filtration $({\cal F}_t)_{t \in [0,T]}$
	satisfying the usual conditions. 
	For all $n \in \{1,\ldots,N\}$, let $
	\xi^{n,m,j}\colon\Omega\to\R^d
	$, 
	$
	m \in \{0, \dots, M-1\}
	$, 
	$ j \in \N $,
	be i.i.d.\ ${\cal F}_0/{\cal B}(\R^d)$-measurable random variables and
	$
	B^{n,m,j} \colon [0,T] \times \Omega \to \R^d 
	$, 
	$
	m \in \{0, \dots, M-1\}
	$, 
	$ 
	j \in \N
	$,
	i.i.d.\ standard
	$({\cal F}_t)_{t \in [0,T]}$-Brownian motions. Let
	$ \bV^{j,{\bf s}}_n\colon\R^\nu\times\R^d \to\R$, $j\in \N$, ${\bf s} \in \R^\varsigma$, $n \in \{0,1,\ldots,N\}$, 
	be functions such that $\bV^{j,{\bf s}}_0(\theta, x) = \varphi(x)$ for all 
	$j\in \N$, ${\bf s} \in \R^\varsigma$, $\theta \in \R^{\nu}$ and $x \in \R^d$.
	For all
	$n \in \{1,\ldots,N\}$,
	$m \in \{0, \dots, M-1\}$ and
	$j \in \N$, let the stochastic process
	$\Y^{n,m,j} \colon \{0,1,\ldots,N\}\times\Omega\to\R^d$
	be given by
	$\Y^{n,m,j}_0 = \xi^{n,m,j}$  and 
	\begin{equation}\label{eq:FormalYapprox}
	\Y^{n,m,j}_{k+1} 
	= 
	H(\tau_{k+1},\tau_{k},\Y^{n,m,j}_k,B^{n,m,j}_{\tau_{k+1}}-B^{n,m,j}_{\tau_{k}}), 
	\quad k\in\{0,1,\ldots,N-1\}.
	\end{equation}  
	For a given function $z \in \cZ$,
	let $\Theta^{z,n}\colon \{0, \dots, M\} \times\Omega\to\R^{\nu}$, $n \in \{1,\dots N\}$, be 
	stochastic processes and $\phi^{z,n,m,{\bf s}} \colon\R^{\nu} \times\Omega\to\R$, 
$n\in\{1,\ldots,N\}$, $m\in \{0, \dots, M-1\}$, ${\bf s}\in\R^{\varsigma}$,
mappings satisfying, for all $\theta \in \R^{\nu}$ and $\omega \in \Omega$,
	\begin{align}
	&\phi^{z,n,m,{\bf s}}(\theta,\omega) 
	= 
	\frac{1}{J_m}\sum_{j=1}^{J_m}
	\bigg[ 
	\bV^{j,{\bf s}}_n\big(\theta,\Y^{n,m,j}_{N-n}(\omega)\big) \nonumber\\
	& - 
	\mathcal{H}_{n-1}\Big(
	\Y^{n,m,j}_{N-n+1}(\omega),
		\bV^{j,{\bf s}}_{n-1}\big(\Theta^{z,n-1}_{M}(\omega),\Y^{n,m,j}_{N-n+1}(\omega)\bigr),\\
		& \quad \quad \quad
	\nabla_x \bV^{j,{\bf s}}_{n-1} \big(\Theta^{z,n-1}_{M}(\omega),
	\Y^{n,m,j}_{N-n+1}(\omega)\big),
	z_{t_n}\big(\Y^{n,m,j}_{N-n+1}(\omega)\big) -
	z_{t_{n-1}}\big(\Y^{n,m,j}_{N-n+1}(\omega)\big)\Big)
	\bigg]^2. \nonumber
	\end{align}
Moreover, assume that for all $n\in\{1,\ldots,N\}$, $m\in \{0, \dots, M-1\}$, 
and ${\bf s}\in\R^{\varsigma}$, $\Phi^{z,n,m,{\bf s}}\colon\R^{\nu}\times\Omega\to\R^{\nu}$
satisfies for all 
$\omega\in\Omega$ and
$\theta\in\{\eta\in\R^{\nu}\colon \phi^{z,n,m,{\bf s}}(\cdot,\omega)
\colon\R^{\nu}\to\R~\text{is differentiable at}~\eta\}$,
	\begin{align}
	\Phi^{z,n,m,{\bf s}}(\theta,\omega) = \nabla_{\theta}\phi^{z,n,m,{\bf s}}(\theta,\omega).
	\end{align}
For $n \in \{1, \dots, N\}$ and $m \in \{0, \dots, M-1\}$, let
$\S^{n}\colon\R^{\varsigma}\times\R^{\nu}\times(\R^d)^{\{0,1,\ldots,N\}\times\N}$ $\to\R^{\varsigma}$, 
	$\psi^{n}_m\colon\R^{\varrho}\to\R^{\nu}$ 
	and $\Psi^{n}_m \colon\R^{\varrho}\times\R^{\nu}\to\R^{\varrho}$
be functions and 
$\bS^{z,n}\colon \{0, \dots, M\} \times\Omega\to\R^{\varsigma}$, 
$\Xi^{z,n}\colon \{0, \dots, M\} \times\Omega\to\R^{\varrho}$ 
stochastic processes such that 
\begin{equation}\label{eq:general_batch_normalization} 
	\bS^{z,n}_{m+1} = \S^{n}\bigl(\bS^{z,n}_m, \Theta^{z,n}_{m}, 
	(\Y_k^{n,m,i})_{(k,i)\in\{0,1,\ldots,N\}\times\N}\bigr),
	\end{equation}
	\begin{equation}
	\Xi^{z,n}_{m+1} = \Psi^{n}_{m}(\Xi^{z,n}_{m},\Phi^{z,n,m,\bS^{z,n}_{m+1}}(\Theta^{z,n}_m))
	\quad
	\text{and}
	\quad
	\Theta^{z,n}_{m+1} = \Theta^{z,n}_{m} - \psi^{n}_{m}(\Xi^{z,n}_{m+1}) 
	\label{eq:general_gradient_step}
	\end{equation}
	for all $n \in \{1, \dots, N\}$ and $m \in \{0, \dots, M-1\}$.
\end{algo}
In the setting of Framework~\ref{def:general_algorithm}, we think for 
sufficiently large $N,M \in \N$, of $\mathbb{V}^{1,\mathbb{S}_M^{z,n}}_n(\Theta^{z,n}_M,x)\colon \Omega \to \R$ as a suitable approximation 
\begin{equation}\label{eq:V-approx-gen-frame}
\mathbb{V}^{1,\mathbb{S}^{z,n}_M}_n(\Theta^{z,n}_M,x) \approx  \E\big[X_{t_n}(x)\!~|\!~Z=z\big],
\end{equation}
where 
$
X \colon [0,T]\times\R^d\times\Omega\to\R
$ 
is a solution of SPDE \eqref{eq:SPDE}.

\section{Examples}
\label{sec:examples}
In this section we illustrate our method on four example SPDEs:
a stochastic heat equation with additive noise,
a stochastic heat equation with multiplicative noise, 
a stochastic Black--Scholes equation with multiplicative noise 
and a Zakai equation.
%
%
In all four cases we use the general  approximation algorithm of Framework~\ref{def:general_algorithm}
in conjunction with the Adam optimizer (see Kingma \& Ba~\cite{KingmaBa2015})
with mini-batches of size 64 in each iteration step 
(see Framework~\ref{frame:adam} below for a detailed description).
We use $N$ fully-connected feedforward neural networks to represent 
$
\bV^{j,{\bf s}}_n(\theta,x) 
$
for $n \in \{ 1, 2, \dots, N\} $,
$ j \in \{ 1, 2, \dots, 64 \} $,
${\bf s} \in \R^{\varsigma}$,
$ \theta \in \R^{ \nu } $ and $x \in \R^d$. 
Each of the neural networks consists of 4 layers
(a $d$-dimensional input layer, two $d+50$-dimensional 
hidden layers and a 1-dimensional output layer).
%
We use $\tanh$ (in Subsections \ref{subsec:stoch_heat}--\ref{subsec:BS}) and GELU 
(in Subsection \ref{subsec:Zakai}) as activation functions in the hidden layers and apply
batch normalization (see Ioffe \& Szegedy~\cite{IoffeSzegedy2015}) just before the first affine transformation (batch normalization for the input) as well as just before every 
activation 
(batch normalization for the hidden layers just before activation).
All parameters of the network are initialized according to a normal or a uniform distribution. 
In Subsections~\ref{subsec:stoch_heat}--\ref{subsec:BS} the mapping $H\colon [0,T]^2\times\R^d\times\R^d\to\R^d$
can be chosen so that the approximation of the auxiliary process \eqref{eq:SDEY} is exact, that is, the processes
$\cY^{n,m,j}$, $n \in \{1, \dots, N\}$, $m \in \{0, \dots, M-1\}$, $j \in \N$, 
from \eqref{eq:FormalYapprox} satisfy 
\[
\cY^{n,m,j}_k = Y^{n,m,j}_{\tau_k} \quad \mbox{for all } k \in \{0,\dots, N\},
\]
where $Y^{n,m,j} \colon [0,T] \times \Omega \to \R$, are strong solutions of 
\be \label{YnmjSDE}
Y^{n,m,j}_t = \xi^{n,m,j} + \int_0^t \mu(Y^{n,m,j}_s ) ds + \int_0^t \sigma(Y^{n,m,j}_s)
dB^{n,m,j}_s, \quad t \in [0,T],
\ee
$n \in \{1, \dots, N\}$, $m \in \{0, \dots, M-1\}$, $j \in \N$. In Subsection~\ref{subsec:Zakai}
we use the Euler--Maruyama scheme to approximate \eqref{YnmjSDE}
(see, e.g., Kloeden \& Platen \cite{KloedenPlaten1992} or Maruyama \cite{Maruyama1955}).

In Subsections~\ref{subsec:stoch_heat} and \ref{subsec:Zakai}, we choose the approximation mappings 
${\cal H}_n \colon \R^d \times \R \times \R^d \times \R^{\zeta} \to \R$, $n \in \{0, \dots, N-1\}$, 
equal to the first order approximations \eqref{1storder}, while in Subsections~\ref{subsec:const-coeff}--\ref{subsec:BS} we use Milstein-type second order approximations
(see, e.g., Kloeden \& Platen \cite[Section~10.3]{KloedenPlaten1992}).

The numerical experiments presented below 
were performed in {\sc Python} using {\sc TensorFlow} on a 
NVIDIA GeForce RTX 4090 GPU. The underlying system 
was an AMD Ryzen~9 7950X CPU with 64 GB DDR5 memory 
running Tensorflow~2.1 on Fedora~40.
The {\sc Python} source codes can be found at \url{https://github.com/seb-becker/deep\_spde}.


\begin{algo}
	\label{frame:adam}
	Assume Framework~\ref{def:general_algorithm} with
	$\nu= (d+50)(d+1)+ (d+50)(d+51)+(d+51)$, $\varepsilon\in (0,\infty)$, 
	$\beta_1 = \tfrac{9}{10}$, 
	$\beta_2 = \tfrac{999}{1000}$, 
	$(\gamma_m)_{m=0}^{M-1}\subseteq (0,\infty)$.
	Let $\operatorname{Pow}_r \colon \R^{\nu}\to\R^{\nu}$ for $r\in (0,\infty)$,
	be given by $\operatorname{Pow}_r(x) = (|x_1|^r,|x_2|^r,\ldots,|x_{\nu}|^r)$.
	Consider functions
	$
	\varphi\colon \R^d \to \R
	$,
	$ 
	\mu=(\mu_1,\mu_2,\dots,\mu_d) \colon \R^d \to \R^d
	$, 
	and 
	$
	\sigma \colon  \R^d \to \R^{ d \times d }.
	$
	Let 
	$ 
	Z=(Z_t(x,\omega))_{(t,x,\omega)\in[0,T]\times\R^d\times \Omega} \colon [0,T] \times \R^d \times \Omega \to \R^\delta
	$ and \linebreak
	$
	X=(X_t(x,\omega))_{(t,x,\omega)\in [0,T]\times \R^d \times \Omega}\colon [0,T]\times\R^d\times\Omega\to\R
	$
	be random fields such for every
	$
	x\in\R^d,
	$ 
	$
	(Z_t(x))_{t\in [0,T]} \colon [0,T] \times \Omega \to \R^\delta
	$ 
	is an 
	$
	(\mathcal F_t)_{t\in [0,T]} $-It\^o process, for all 
	$
	t\in [0,T]
	$ and
	$
	x\in\R^d,
	$ 
	$
	X_t(x)\colon\Omega\to\R
	$ 
	is 
	$\mathcal F_t$/$\B(\R)$-measurable, 
%
	for all $\omega\in\Omega$,
	$
	(X_t(x,\omega))_{(t,x) \in [0,T]\times\R^d}$ is in $C^{0,2}([0,T]\times\R^d,\R)$
with at most polynomially growing partial derivatives of order 0, 1 and 2
with respect to the $x$-variables, and for all
	$t \in [0,T]$ and $x\in\R^d$, one has 
	\begin{align*}
	X_{ t }( x )
	&  =
	\varphi( x ) 
	+
	\int_{ 0 }^{ t }
	f \big( 
	x, X_s(x), \nabla X_s( x ) 
	\big)
	\, ds
	+
	\int_{ 0 }^{ t }
	\big\langle b\big( x, X_s( x ), \nabla X_s( x ) \big), dZ_s(x) \big\rangle_{\R^\delta}
\nonumber	\\
	& \quad
	+
	\int_{ 0 }^{ t }
	\Big(
	\tfrac{ 1 }{ 2 }
	\Tr \big( 
	\sigma( x ) \sigma(x )^T
	\He X_s ( x )
	\big)
	+
	\big\langle \mu( x ), \nabla X_s ( x ) \big\rangle_{ \R^d }
	\Big)
	\, ds \; \; \mbox{$\P$-a.s.}
	\end{align*}
	Let $ \varrho = 2 \nu $, and assume $ t_i = \tfrac{iT}{N} $ for $i \in \{0, \dots, N\}$
	and $ J_m = 64 $ for all $m \in \{0, \dots, M\}$. Moreover, let 
	\begin{align}\label{eq:examples_setting_moment_estimation}
	\Psi^n_m ( x , y , \eta ) 
	= 
	(\beta_1 x + (1-\beta_1) \eta, \beta_2 y + (1-\beta_2) \operatorname{Pow}_2(\eta))
	\end{align}
	and 
	\begin{align}\label{eq:examples_setting_adam_grad_update}
	\psi^n_m ( x,y ) = 
	\biggl(
	\Bigl[
	\sqrt{\tfrac{|y_1|}{1-(\beta_2)^m}} + \varepsilon
	\Bigr]^{-1}
	\frac{\gamma_m x_{1}}{1-(\beta_1)^m},
	\ldots, 
	\Bigl[
	\sqrt{\tfrac{|y_{\nu}|}{1-(\beta_2)^m}} + \varepsilon
	\Bigr]^{-1}
	\frac{\gamma_m x_{\nu}}{1-(\beta_1)^m}
	\biggr). 
	\end{align}
	for $x, y, \eta \in \R^{\nu}$.	
\end{algo}
$J_m = 64$ is the batch size, and equations 
\eqref{eq:examples_setting_moment_estimation}--\eqref{eq:examples_setting_adam_grad_update} describe the Adam optimizer (see Kingma \& Ba \cite{KingmaBa2015}.

\subsection{Stochastic heat equation with  additive noise}
\label{subsec:stoch_heat}

Our first example is the stochastic heat equation with additive noise
\begin{equation}\label{eq:ex-heat-add}
	X_t(x)
	= 
	\|x\|_{\R^d}^2
	+ 
	\int_0^t 
	\Delta_x X_s(x) ds
	+ 
	W_t, \quad (t,x) \in [0,T] \times \R^d,
\end{equation}
for a standard Brownian motion $W\colon [0,T]\times\Omega\to\R$. 
It can easily be checked that it admits 
the following closed-form solution 
\begin{equation}\label{le:eq:SPDE-heat-add}
X_t(x) = \|x\|_{\R^d}^2 + 2td + W_t, \quad (t,x) \in [0,T] \times \R^d,
\end{equation}
which allows us to obtain reference values for our numerical experiments. 

To numerically solve \eqref{eq:ex-heat-add}, we used Framework~\ref{frame:adam} with
$T=1$,
$N=5$, 
$\varepsilon=10^{-8}$,
$\varphi(x)=\|x\|_{\R^d}^2$, $f(x,u,v) = 0$, $b(x,u,v) = 1$, 
$\mu(x)=0$, $\sigma(x)v=\sqrt{2}v$, 
$\delta = 1$,
$Z_t(x) = W_t$,
$H(t,s,x,v)= x +\sqrt{2} v$,
$\mathcal{H}_n(x,u,v,w)=u+w$ for $n \in \{0, \dots, N-1\}$, $x,v \in \R^d$, $u,w \in \R$ and
$s,t \in [0,T]$ such that $s < t$.

To approximate the realization $X^z_T(0)$ of $X_T(0)$ along a given path $z$ of $W$,
we set $\xi^{n,m,j}= 0 \in \R^d$ for all $n \in \{1, \dots, N \}$, $m \in \{0, \dots, M-1\}$,
$j \in \N$, $M=10000$ and $\gamma_m = 10^{-1} \mathbbm{1}_{[0,4000]}(m) 
+ 10^{-2} \mathbbm{1}_{(4000,6000]}(m)
+ 10^{-3} \mathbbm{1}_{(6000,8000]}(m)
+ 10^{-4} \mathbbm{1}_{(8000,10000)}(m)$ for $m \in \{0,1,\dots, M-1\}$.
Table \ref{table:heat-add} shows numerical results for five different realizations $z$ of 
$W$ in each of the cases $d \in \{1, 5, 10, 20, 50, 100\}$. 
It reports the results $\mathbb{V}^{1,\mathbb{S}^{z,N}_M}_N(\Theta^{z,N}_M,0)$
of our approximation algorithm, computation times, reference solutions 
$X^z_T(0)$ obtained from \eqref{le:eq:SPDE-heat-add} along the same realizations 
$z$ of $W$, pathwise relative errors 
\begin{equation} \label{relerror}
\frac{
  \bigl|
    \mathbb{V}_N^{1,\mathbb{S}^{z,N}_M}(\Theta^{z,N}_M,0) - X^z_T(0)
  \bigr|
}{
  \bigl| X^z_T(0) \bigr|
}
\end{equation}
and averages of \eqref{relerror} over the five different realizations $z$ of $W$.

\begin{table}
	\begin{center}
		\small
		\begin{tabular}{|c|c|c|c|c|c|}
			\hline
			$d$ & \makecell{Result of\\ the approx.\\algorithm} & \makecell{Runtime\\ in\\ seconds} & \makecell{Reference\\ solution} & \makecell{Pathwise \\ relative\\ error} & \makecell{Avg.\\ relative\\ error}\\
			\hline
			  & 1.618 & 33.83 & 1.596 & 0.0139 & \\ 
			  & 3.144 & 33.84 & 3.176 & 0.0102 & \\ 
			1 & 0.452 & 33.76 & 0.454 & 0.0053 & 0.0010 \\ 
			  & 3.770 & 33.16 & 3.763 & 0.0019 & \\ 
			  & 2.403 & 33.41 & 2.448 & 0.0185 & \\ 
			\hline
			  & 12.141 & 34.00 & 12.083 & 0.0048 & \\ 
			  & 9.008 & 33.59 & 8.909 & 0.0112 & \\ 
			5 & 9.135 & 34.64 & 9.080 & 0.0061 & 0.0064 \\ 
			  & 9.149 & 34.39 & 9.088 & 0.0067 & \\ 
			  & 12.012 & 33.77 & 11.976 & 0.0030 & \\ 
			\hline
			   & 19.181 & 34.46 & 19.220 & 0.0021 & \\ 
			   & 20.905 & 33.97 & 20.844 & 0.0030 & \\ 
			10 & 20.664 & 34.05 & 20.577 & 0.0042 & 0.0024 \\ 
			   & 19.861 & 33.50 & 19.874 & 0.0006 & \\ 
			   & 20.101 & 34.50 & 20.143 & 0.0021 & \\ 
			\hline
			   & 39.969 & 34.42 & 40.020 & 0.0013 & \\ 
			   & 40.372 & 34.38 & 40.091 & 0.0070 & \\ 
			20 & 39.671 & 34.32 & 39.474 & 0.0050 & 0.0044 \\ 
			   & 40.082 & 34.35 & 39.810 & 0.0068 & \\ 
			   & 39.968 & 34.07 & 39.893 & 0.0019 & \\ 
			\hline
			   & 99.945 & 34.75 & 99.935 & 0.0001 & \\ 
			   & 100.426 & 34.46 & 100.096 & 0.0033 & \\ 
			50 & 98.395 & 34.62 & 98.226 & 0.0017 & 0.0014 \\ 
			   & 99.674 & 34.49 & 99.743 & 0.0007 & \\ 
			   & 99.116 & 34.90 & 99.005 & 0.0011 & \\ 
			\hline
			    & 199.039 & 35.80 & 201.098 & 0.0102 & \\ 
			    & 199.385 & 36.31 & 201.236 & 0.0092 & \\ 
			100 & 198.890 & 35.98 & 200.964 & 0.0103 & 0.010 \\ 
			    & 197.945 & 35.84 & 200.078 & 0.0107 & \\ 
			    & 198.870 & 35.67 & 200.698 & 0.0091 & \\ 
			\hline
		\end{tabular}
		\caption{\sl Numerical results for the stochastic heat equation with 
		additive noise \eqref{eq:ex-heat-add}.
		}
		\label{table:heat-add}
	\end{center}
\end{table}

\begin{figure}[h]
\begin{center}
\includegraphics[width=0.6\textwidth]{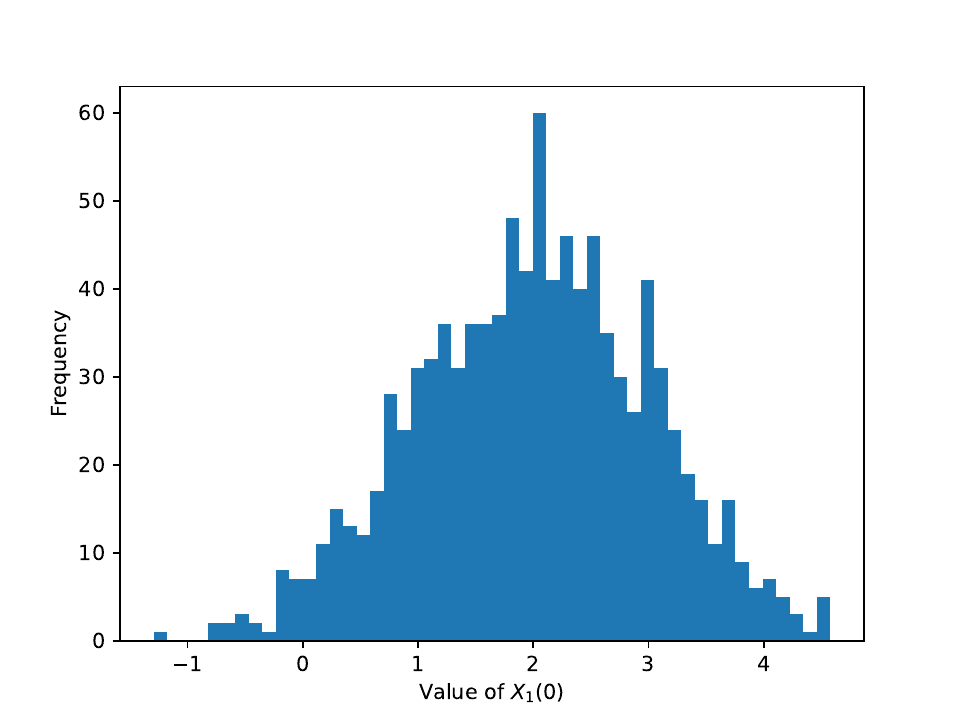}
\caption{\sl Histogram of 1000 approximate realizations of $X_1(0)$ for the stochastic heat equation with additive noise
\eqref{eq:ex-heat-add} in the case $d=1$.}
\label{fig:histo}
\end{center}
\end{figure}

The distribution of $X_T(x)$ can be approximated by running the algorithm on a larger set of 
simulated trajectories of the noise process $Z$. This yields an empirical distribution from which 
functionals such as the mean, variance or standard deviation can be estimated. Figure \ref{fig:histo}
shows a histogram of numerical approximations of $X^z_1(0)$ in the case $d=1$ 
for 1000 simulated noise trajectories $z$ of $Z$. Table \ref{table:histo} shows sample based estimates of 
the mean $\E [X_1(0)]$ and standard deviation $\sqrt{{\rm Var}(X_1(0))}$ together 
with their true values obtained from the explicit solution \eqref{le:eq:SPDE-heat-add}.
It can be seen that the estimates are very accurate, but the computation takes a long time
because $X^z_T(0)$ must be approximated for 1000 realizations $z$ of $Z$.

\begin{table}[ht]
	\begin{center}
		\small
		\begin{tabular}{|c|c|c|c|c|}
			\hline
			 \makecell{Sample\\ mean} & \makecell{True\\ mean} & \makecell{Sample standard\\ deviation} & \makecell{True standard\\ deviation} & \makecell{Runtime}\\
			\hline
			   1.9971 & 2 & 1.0079 & 1 & 7h 39min 16s\\
			\hline
		\end{tabular}
		\caption{\sl Estimates of the mean and standard deviation 
		of $X_1(0)$ compared to their true values for the stochastic heat equation with 
		additive noise \eqref{eq:ex-heat-add} in the case $d=1$.
		}
		\label{table:histo}
	\end{center}
\end{table}

To obtain approximations of the function $x \mapsto X^z_T(x)$ for a given realization
$z$ of $W$, one can choose the initial conditions $\xi^{n,m,j}$, $n \in \{1, \dots, N \}$, 
$m \in \{0, \dots, M-1\}$, $j \in \N$, 
randomly instead of deterministically. For the results of Table \ref{table:heat-add-domain} 
we chose them i.i.d.\ uniformly distributed over $[-1,1]^d$ and set $M = 20000$ together with
$\gamma_m = 10^{-1} \mathbbm{1}_{[0,5000]}(m) 
+ 10^{-2} \mathbbm{1}_{(5000,10000]}(m)
+ 10^{-3} \mathbbm{1}_{(10000,15000]}(m)
+ 10^{-4} \mathbbm{1}_{(15000,20000)}(m)$
for $m \in \{0,1, \dots, M-1\}$. Table \ref{table:heat-add-domain} shows for
five different realizations $z$ of $W$ in each of the cases
$d \in \{1,5,10,20,50,100\}$, Monte Carlo estimates of the $L^1$-error
\begin{equation} \label{L1error}
\frac{1}{2^d} \int_{[-1,1]^d} 
\abs{\mathbb{V}_N^{1,\mathbb{S}^{z,N}_M}\!\!(\Theta^{z,N}_M,x)-X^z_T(x)} dx,
\end{equation}
the $L^1$-norm of the relative errors
\begin{equation}
\frac{1}{2^d} \int_{[-1,1]^d} \frac{\abs{\mathbb{V}_N^{1,\mathbb{S}^{z,N}_M}\!\!
(\Theta^{z,N}_M,x)-X^z_T(x)}}{|X^z_T(x)|} \, dx,
\end{equation} 
the $L^2$-error
\begin{equation} \label{L2error}
\brak{\frac{1}{2^d} \int_{[-1,1]^d} 
{\abs{\mathbb{V}_N^{1,\mathbb{S}^{z,N}_M}\!\!(\Theta^{z,N}_M,x)-X^z_T(x)}^2}\,dx }^{1/2},
\end{equation}
the $L^2$-norm of the relative errors
\begin{equation} \label{L2norm}
\brak{\frac{1}{2^d} \int_{[-1,1]^d} \brak{\frac{\mathbb{V}_N^{1,\mathbb{S}^{z,N}_M}\!\!
(\Theta^{z,N}_M,x)-X^z_T(x)}{X^z_T(x)}}^2 \, dx}^{1/2}
\end{equation}
and average computation times.

\begin{table}
	\begin{center}
		\small
		\begin{tabular}{|c|c|c|c|c|c|}
			\hline
			$ d $ & $L^1$-error & \makecell{$L^1$-norm of\\ rel. errors}
			& $L^2$-error & \makecell{$L^2$-norm of\\ rel. errors} & \makecell{Avg.\ runtime\\ 
			in seconds}\\
			\hline
			  & 0.007066 & 0.002856 & 0.009375 & 0.003676 & \\
			  & 0.013098 & 0.012499 & 0.018331 & 0.016024 & \\
			1 & 0.011375 & 0.006832 & 0.013698 & 0.008084 & 74.2\\
			  & 0.012296 & 0.004263 & 0.015563 & 0.005245 & \\
			  & 0.016615 & 0.012314 & 0.021642 & 0.015112 & \\
			\hline
			  & 0.035044 & 0.002658 & 0.044177 & 0.003328 & \\
			  & 0.043454 & 0.003699 & 0.054438 & 0.004589 & \\
			5 & 0.036685 & 0.003226 & 0.045460 & 0.003985 & 75.4\\
			  & 0.059150 & 0.005733 & 0.067932 & 0.006594 & \\
			  & 0.043178 & 0.004043 & 0.053583 & 0.004994 & \\
			\hline
			   & 0.069180 & 0.002993 & 0.085204 & 0.003691 & \\
			   & 0.083577 & 0.003503 & 0.102072 & 0.004277 & \\
			10 & 0.057488 & 0.002704 & 0.072292 & 0.003384 & 75.6 \\
			   & 0.076894 & 0.003504 & 0.092920 & 0.004238 & \\
			   & 0.055275 & 0.002293 & 0.069274 & 0.002866 & \\
			\hline
			   & 0.107869 & 0.002222 & 0.134527 & 0.002773 & \\
			   & 0.120857 & 0.002499 & 0.149345 & 0.003090 & \\
			20 & 0.120547 & 0.002683 & 0.149282 & 0.003323 & 76.5 \\
			   & 0.098655 & 0.002127 & 0.123768 & 0.002664 & \\
			   & 0.108540 & 0.002309 & 0.136028 & 0.002890 & \\
			\hline
			   & 0.607528 & 0.005234 & 0.751561 & 0.006477 & \\
			   & 0.611149 & 0.005211 & 0.754792 & 0.006439 & \\
			50 & 0.728106 & 0.006320 & 0.898680 & 0.007812 & 85.5 \\
			   & 0.607731 & 0.005258 & 0.753373 & 0.006518 & \\
			   & 0.591882 & 0.005073 & 0.739039 & 0.006325 & \\
			\hline
			    & 1.169286 & 0.004994 & 1.458186 & 0.006231 & \\
			    & 1.230534 & 0.005276 & 1.531159 & 0.006573 & \\
			100 & 1.050782 & 0.004490 & 1.311571 & 0.005604 & 105.3 \\
			    & 1.043714 & 0.004494 & 1.310187 & 0.005635 & \\
			    & 1.098031 & 0.004712 & 1.368130 & 0.005873 & \\
			\hline
		\end{tabular}
		\caption{\sl $L^1$- and $L^2$-errors for the stochastic heat equation with additive noise
			 \eqref{eq:ex-heat-add}.
		}
	\label{table:heat-add-domain}
	\end{center}
\end{table}


\subsection{Stochastic heat equation with multiplicative noise}
\label{subsec:const-coeff}

Next, we consider the stochastic heat equation with multiplicative noise
\begin{equation}\label{eq:ex-heat-multi}
	X_t(x)
	= 
	\|x\|_{\R^d}^2
	+ 
	\int_0^t 
	\Delta_x X_s(x) ds
	+ \int_0^t X_s(x)
	dW_s, \quad (t,x) \in [0,T] \times \R^d,
\end{equation}
for a standard Brownian motion $W\colon [0,T]\times\Omega\to\R$. An application of 
It\^{o}'s formula shows that it admits the closed-form solution
\begin{equation} \label{le:eq:SPDE-heat-multi}
	X_t(x) = \big( \|x\|_{\R^d}^2 +2 td \big) e^{W_t - \frac{1}{2} t}, 
	\quad (t,x) \in [0,T] \times \R^d,
	\end{equation} which we use to produce reference values 
for our numerical results.

To obtain a numerical solution to \eqref{eq:ex-heat-multi}, 
we adopted Framework~\ref{frame:adam} with 
$T=0.5$, $N=25$, $\varepsilon=10^{-8}$, 
$\varphi(x)=\|x\|_{\R^d}^2$, $f(x,u,v) = 0$, $b(x,u,v) = u$,
$\mu(x)=0$, $\sigma(x)v=\sqrt{2}v$, 
$\delta = 1$, $Z_t(x) = W_t$,
$H(t,s,x,v)= x +\sqrt{2} v$, $\mathcal{H}_n(x,u,v,w)=u\big(1+w+\frac{1}{2}w^2-\frac{1}{2}(t_{n+1}-t_n)\big)$ 
for $n \in \{0, \dots, N-1\}$, $x,v \in \R^d$, $u,w \in \R$ and $s, t \in [0,T]$ such that $s < t$.

First, we set $\xi^{n,m,j}= 0 \in \R^d$ for all $n \in \{1, \dots, N \}$, $m \in \{0, \dots, M-1\}$,
$j \in \N$ and $M=12000$ together with $\gamma_m = 10^{-1} \mathbbm{1}_{[0,5000]}(m) 
+ 10^{-2} \mathbbm{1}_{(5000,7000]}(m)
+ 10^{-3} \mathbbm{1}_{(7000,10000]}(m)
+ 10^{-4} \mathbbm{1}_{(10000,12000)}(m)$ for $m \in \{0,1,\dots, M-1\}$.
Table \ref{table:heat-mult} shows results 
$\mathbb{V}^{1,\mathbb{S}^{z,N}_M}_N(\Theta^{z,N}_M,0)$
of our approximation algorithm, computation times, reference solutions 
$X^z_T(0)$ obtained from \eqref{le:eq:SPDE-heat-multi}, pathwise relative errors 
\be \label{relerror-multi}
\frac{
  \bigl|
    \mathbb{V}_N^{1,\mathbb{S}^{z,N}_M}(\Theta^{z,N}_M,0) - X^z_T(0)
  \bigr|
}{
  \bigl|
    X^z_T(0)
  \bigr|
}
\ee
for five different realizations $z$ of $W$ in each of the cases
$d \in \{1, 5, 10, 20, 50, 100\}$ as well as averages of
\eqref{relerror-multi} over the five different realizations $z$ of $W$.

\begin{table}
	\begin{center}
		\small
			\begin{tabular}{|c|c|c|c|c|c|}
				\hline
				\makecell{$d$} & \makecell{Result of\\ the approx.\\algorithm} & \makecell{Runtime\\ in\\ seconds} & \makecell{Reference\\ solution} & \makecell{Pathwise\\ relative\\ error} & \makecell{Avg.\\ relative\\ error}\\
				\hline
				  & 1.193 & 225.33 & 1.180 & 0.0110 & \\ 
				  & 0.873 & 225.91 & 0.884 & 0.0119 & \\ 
				1 & 0.647 & 225.31 & 0.623 & 0.0372 & 0.0183 \\ 
				  & 0.598 & 224.00 & 0.597 & 0.0015 & \\ 
				  & 1.153 & 225.97 & 1.189 & 0.0299 & \\ 
				\hline
				  & 1.034 & 230.27 & 1.022 & 0.0120 & \\ 
				  & 1.213 & 229.59 & 1.192 & 0.0173 & \\ 
				5 & 4.673 & 228.43 & 4.629 & 0.0096 & 0.0117 \\ 
				  & 11.647 & 230.20 & 11.422 & 0.0197 & \\ 
				  & 4.709 & 230.56 & 4.708 & 0.0001 & \\ 
				\hline
				   & 7.383 & 231.44 & 7.343 & 0.0054 & \\ 
				   & 3.129 & 230.08 & 3.092 & 0.0117 & \\ 
				10 & 13.812 & 231.98 & 13.757 & 0.0040 & 0.0101 \\ 
				   & 11.065 & 231.24 & 10.902 & 0.0150 & \\ 
				   & 7.110 & 230.12 & 7.010 & 0.0142 & \\ 
				\hline
				   & 29.450 & 229.46 & 28.794 & 0.0228 & \\ 
				   & 12.045 & 229.01 & 12.002 & 0.0036 & \\ 
				20 & 5.381 & 228.11 & 5.383 & 0.0003 & 0.0064 \\ 
				   & 16.594 & 228.39 & 16.519 & 0.0045 & \\ 
				   & 10.935 & 231.57 & 10.928 & 0.0007 & \\ 
				\hline
				   & 25.865 & 232.01 & 25.697 & 0.0065 & \\ 
				   & 189.230 & 230.36 & 188.450 & 0.0041 & \\ 
				50 & 51.230 & 232.80 & 51.017 & 0.0042 & 0.0045 \\ 
				   & 169.914 & 232.88 & 168.914 & 0.0059 & \\ 
				   & 55.788 & 232.55 & 55.677 & 0.0020 & \\ 
				\hline
				    & 42.704 & 254.85 & 42.779 & 0.0018 & \\ 
				    & 114.128 & 254.80 & 109.905 & 0.0384 & \\ 
				100 & 58.949 & 255.52 & 57.983 & 0.0167 & 0.0134 \\ 
				    & 90.850 & 253.83 & 91.201 & 0.0038 & \\ 
				    & 246.615 & 255.38 & 248.238 & 0.0065 & \\ 
				\hline
			\end{tabular}
		\caption{\sl Numerical results for the stochastic heat equation with multiplicative noise \eqref{eq:ex-heat-multi}.
		}
		\label{table:heat-mult}
	\end{center}
\end{table}

To derive approximations of the function $x \mapsto X^z_T(x)$ for given realizations
$z$ of $W$, we chose $\xi^{n,m,j}$, $n \in \{1, \dots, N \}$, $m \in \{0, \dots, M-1\}$,
$j \in \N$, i.i.d.\ uniformly distributed over $[-1,1]^d$ and set 
$M=20000$ together with 
$\gamma_m = 10^{-1} \mathbbm{1}_{[0,5000]}(m) 
+ 10^{-2} \mathbbm{1}_{(5000,10000]}(m)
+ 10^{-3} \mathbbm{1}_{(10000,15000]}(m)
+ 10^{-4} \mathbbm{1}_{(15000,20000)}(m)$, $m \in \{0,1,\dots, M-1\}$.
Table \ref{table:heat-mult-domain} shows Monte Carlo estimates of 
\eqref{L1error}--\eqref{L2norm} for five different realizations $z$ of 
$W$ and average runtimes in each of the cases $d \in \{1,5,10,20,50,100\}$.

\begin{table}
	\begin{center}
		\small
		\begin{tabular}{|c|c|c|c|c|c|}
			\hline
			$ d $ & $L^1$-error & \makecell{$L^1$-norm of\\ rel. errors} 
			& $L^2$-error & \makecell{$L^2$-norm of\\ rel. errors} & 
			\makecell{Avg.\ runtime\\ in seconds}\\
			\hline
			  & 0.099791 & 0.030543 & 0.117261 & 0.034766 & \\
			  & 0.007980 & 0.009697 & 0.009616 & 0.011233 & \\
			1 & 0.015422 & 0.029434 & 0.019377 & 0.035182 & 390.0 \\
			  & 0.048228 & 0.028835 & 0.052479 & 0.031079 & \\
			  & 0.018211 & 0.022474 & 0.022648 & 0.026368 & \\
			\hline
			  & 0.038465 & 0.007646 & 0.047532 & 0.009333 & \\
			  & 0.073100 & 0.006760 & 0.091737 & 0.008372 & \\
			5 & 0.016873 & 0.006236 & 0.021054 & 0.007728 & 391.5 \\
			  & 0.122509 & 0.018564 & 0.133084 & 0.020035 & \\
			  & 0.055346 & 0.009671 & 0.070375 & 0.012264 & \\
			\hline
			   & 0.076959 & 0.007700 & 0.091598 & 0.009067 & \\
			   & 0.207491 & 0.010611 & 0.236149 & 0.012077 & \\
			10 & 0.171416 & 0.008851 & 0.197242 & 0.010286 & 394.8 \\
			   & 0.075860 & 0.005587 & 0.092655 & 0.006834 & \\
			   & 0.047575 & 0.013862 & 0.052554 & 0.015076 & \\
			\hline
			   & 0.497157 & 0.022897 & 0.506898 & 0.023375 & \\
			   & 0.273638 & 0.014125 & 0.288383 & 0.014911 & \\
			20 & 0.310580 & 0.008673 & 0.347467 & 0.009710 & 398.0 \\
			   & 0.100438 & 0.007693 & 0.115465 & 0.008766 & \\
			   & 0.116730 & 0.004485 & 0.146831 & 0.005599 & \\
			\hline
			   & 0.233946 & 0.006809 & 0.277234 & 0.008020 & \\
			   & 0.762201 & 0.007644 & 0.872832 & 0.008790 & \\
			50 & 0.676821 & 0.006375 & 0.790546 & 0.007482 & 442.4 \\
			   & 0.259814 & 0.008870 & 0.287583 & 0.009853 & \\
			   & 0.411384 & 0.010830 & 0.444561 & 0.011747 & \\
			\hline
			    & 0.473851 & 0.004572 & 0.582707 & 0.005647 & \\
			    & 0.607502 & 0.004585 & 0.756786 & 0.005717 & \\
			100 & 0.332232 & 0.004173 & 0.419929 & 0.005253 & 614.7 \\
			    & 0.174889 & 0.004983 & 0.212811 & 0.006092 & \\
			    & 1.629214 & 0.007585 & 1.893700 & 0.008768 & \\
			\hline
		\end{tabular}
		\caption{\sl $L^1$- and $L^2$-errors for the stochastic heat equation with 
		multiplicative noise \eqref{eq:ex-heat-multi}.
		}
		\label{table:heat-mult-domain}
	\end{center}
\end{table}

\subsection{Stochastic Black--Scholes equation with multiplicative noise}
\label{subsec:BS}

In this subsection we apply our method to the 
stochastic Black--Scholes-type equation with multiplicative noise
\begin{equation}
\label{eq:BS}
X_t(x)
= 
\varphi(x) 
+ 
\int_0^{t} 
\bigg(
\tfrac12 {\textstyle\sum\limits_{i=1}^d} \sigma_i^2 x_i^2 \tfrac{\partial^2}{\partial x^2_i}X_s(x) 
+ 
{\textstyle\sum\limits_{i=1}^d} \mu_i x_i \tfrac{\partial}{\partial x_i}X_s(x) 
\bigg) ds
+ 
\int_0^t X_s(x) dW_s, \quad (t,x) \in [0,T] \times \R^d,
\end{equation}
for $\varphi(x) = \exp(-rT) \max \big\{\max_{i\in\{1,2,\dots,d\}} x_i -100,0\big\}$ with
$T = 0.5$, $r = 0.02$, $\mu_1 = \nicefrac{\sin(d) + 1}{d}$, 
$\mu_2 = \nicefrac{\sin(2d) + 1}{d}, \dots, \mu_d 
= \nicefrac{\sin(d^2) + 1}{d}$ and $\sigma_1 = \nicefrac{1}{4d}$,
$\sigma_2= \nicefrac{2}{4d}, \dots , \sigma_d=\nicefrac{d}{4d}$. Note that
the function $p \colon [0,T] \times \R^d \to \R$ given by 
$p(0,x)=\varphi(x)$, $x \in \R^d$, and
\[
\begin{split}
p(t,x) = &
\tfrac{1}{(2\pi t)^{\nicefrac{d}{2}}} 
\int_{\R} \int_\R \dots \int_\R 
\Bigg[\exp\!\bigg(-\frac{\sum_{i=1}^d y_i^2}{2t}\bigg)\\
& \varphi\bigg(x_1 e^{\sigma_1 y_1 + [\mu_1- \sigma_1^2/2]t},\ldots,
x_d e^{\sigma_dy_d + [\mu_d - \sigma_d^2/2]t}\bigg)\Bigg]\,dy_1 \,dy_2 \,\dots dy_d,
\; (t,x) \in (0,T] \times \R^d, 
\end{split}
\]
can be written as 
\be \label{pE}
p(t,x) = 
\E \edg{ \varphi \brak{x_1 e^{\sigma_1 B^1_t + [\mu_1 - \sigma^2_1/2]t}, 
\dots, x_d e^{\sigma_d B^d_t + [\mu_d - \sigma^2_d/2]t}}}, \quad 
(t,x) \in [0,T] \times \R^d,
\ee
for a $d$-dimensional standard Brownian motion $B \colon [0,T] \times \Omega \to \R^d$. 
Since $\varphi \colon \R^d \to \R$ is 
continuous and of linear growth, it can be seen from \eqref{pE} that
$p$ is continuous on $[0, T] \times \R^d$, and one obtains from
the Feynman--Kac formula that 
\[
		\label{eq:BS-PDE}
		\frac{\partial}{\partial t} p (t,x)= 
		\frac12 \sum_{i=1}^d \sigma_i^2 x_i^2 \frac{\partial^2}{\partial x^2_i} p(t,x) 
		+ 
		\sum_{i=1}^d \mu_i x_i \frac{\partial}{\partial x_i}p (t,x)
\quad \mbox{for all } (t,x) \in (0,T] \times \R^d,
		\]
which, together with It\^{o}'s formula, shows that
\begin{equation} \label{BSsol}
X_t(x) = p(t,x) e^{W_t - \frac{1}{2} t}, \quad t \in [0,\infty), \, x \in \R^d,
\end{equation}
solves \eqref{eq:BS}. 
So we can use \eqref{BSsol} to compute reference solutions.

To approximate the solution of equation \eqref{eq:BS} 
with our method, we used Framework~\ref{frame:adam} with
$T=0.5$, $N=20$, $\varepsilon=10^{-8}$,
$\varphi(x) = \exp(-rT)\max\!\big\{\max_{i\in\{1,2,\dots,d\}} x_i -100,0\big\}$,
$
f(x,u,v) = 0
$,
$
b(x,u,v) = u
$,
$\langle \mu(x), v \rangle_{\R^d} = \sum_{i=1}^d \mu_i x_i v_i$,
$\sigma(x)v= (\sigma_1 x_1 v_1, \dots, \sigma_d x_d v_d)$,
$\delta = 1$,
$Z_t(x) = W_t$,
\begin{equation}
H(t,s,x,v)= 
\left( x_1 \exp \left([\mu_1 - \sigma_1^2/2](t-s) + \sigma_1 v_1 \right), \dots, 
x_d \exp \left([\mu_d - \sigma_d^2/2](t-s) + \sigma_d v_d \right) \right)
\end{equation} 
and $\mathcal{H}_n(x,u,v,w)= u(1 + w + \frac{1}{2} w^2 - \frac{1}{2} (t_{n+1} - t_n))$ 
for $n \in \{0, \dots, N-1\}$, $x,v \in \R^d$, $u,w \in \R$ and $s,t \in [0,T]$ such that $s < t$.

We first approximated $X^z_T(x)$ for different realizations $z$ of $W$ and 
$x = (100, \dots, 100) \in \R^d$ for $d \in \{1, 5, 10, 20, 50, 100\}$. To this end, we
set $\xi^{n,m,j}= (100, \dots, 100) \in \R^d$ for all $n \in \{1, \dots, N \}$, $m \in \{0, \dots, M-1\}$,
$j \in \N$ and $M=10000$ together with $\gamma_m = 10^{-1} \mathbbm{1}_{[0,4000]}(m) 
+ 10^{-2} \mathbbm{1}_{(4000,6000]}(m)
+ 10^{-3} \mathbbm{1}_{(6000,8000]}(m)
+ 10^{-4} \mathbbm{1}_{(8000,10000)}(m)$ for $m \in \{0,1,\dots, M-1\}$.
Table \ref{table:black-mult}  shows results \linebreak
$\mathbb{V}^{1,\mathbb{S}^{z,N}_M}_N(\Theta^{z,N}_M,100, \dots, 100)$
of our approximation algorithm, computation times, reference solutions 
$X^z_T(100, \dots, 100)$ obtained from \eqref{BSsol} and pathwise relative errors 
\begin{equation} \label{relerror-BS}
\frac{
  \bigl|
    \mathbb{V}_N^{1,\mathbb{S}^{z,N}_M}(\Theta^{z,N}_M,100, \dots, 100) - X^z_T(100, \dots, 100)
  \bigr|
}{
  \bigl| 
    X^z_T(100, \dots, 100)
  \bigr|
}
\end{equation}
for five different realizations $z$ of $W$ in each of the cases
$d \in \{1, 5, 10, 20, 50, 100\}$. The last column shows averages of the relative errors
\eqref{relerror-BS} over the five different realizations $z$ of $W$.

\begin{table}[ht]
	\begin{center}
		\small
		\begin{tabular}{|c|c|c|c|c|c|}
			\hline
			\makecell{$d$} & \makecell{Result of\\ the approx.\\ algorithm} & \makecell{Runtime\\ in\\ seconds} & \makecell{Reference\\ solution} & \makecell{Pathwise\\ relative\\ error} & \makecell{Avg.\\ relative\\ error}\\
			\hline
			  & 141.589 & 153.82 & 142.785 & 0.0084 & \\ 
			  & 140.028 & 152.55 & 140.584 & 0.0040 & \\ 
			1 & 90.429 & 151.82 & 89.448 & 0.0110 & 0.0075 \\ 
			  & 53.750 & 150.38 & 54.392 & 0.0118 & \\ 
			  & 208.328 & 152.76 & 208.847 & 0.0025 & \\ 
			\hline
			  & 6.701 & 154.71 & 6.768 & 0.0099 & \\ 
			  & 45.195 & 153.90 & 44.623 & 0.0128 & \\ 
			5 & 27.949 & 154.36 & 28.441 & 0.0173 & 0.0149 \\ 
			  & 38.273 & 154.54 & 38.717 & 0.0115 & \\ 
			  & 15.042 & 153.70 & 14.707 & 0.0228 & \\ 
			\hline
			   & 15.843 & 155.67 & 15.898 & 0.0035 & \\ 
			   & 33.740 & 155.05 & 33.873 & 0.0039 & \\ 
			10 & 13.928 & 155.95 & 14.061 & 0.0095 & 0.0055 \\ 
			   & 47.863 & 156.45 & 47.802 & 0.0013 & \\ 
			   & 9.846 & 156.73 & 9.754 & 0.0094 & \\ 
			\hline
			   & 18.951 & 156.22 & 18.253 & 0.0382 & \\ 
			   & 17.769 & 158.00 & 17.684 & 0.0048 & \\ 
			20 & 11.657 & 157.52 & 11.460 & 0.0171 & 0.0204 \\ 
			   & 12.246 & 155.71 & 11.984 & 0.0218 & \\ 
			   & 12.806 & 156.04 & 12.557 & 0.0199 & \\ 
			\hline
			   & 30.342 & 156.66 & 28.959 & 0.0477 & \\ 
			   & 26.603 & 158.03 & 25.724 & 0.0342 & \\ 
			50 & 28.645 & 155.12 & 27.649 & 0.0360 & 0.0390 \\ 
			   & 19.954 & 156.62 & 19.298 & 0.0340 & \\ 
			   & 38.748 & 157.69 & 37.149 & 0.0430 & \\ 
			\hline
			    & 85.939 & 161.16 & 82.265 & 0.0447 & \\ 
			    & 14.276 & 163.86 & 13.690 & 0.0428 & \\ 
			100 & 55.103 & 161.44 & 53.146 & 0.0368 & 0.0422 \\ 
			    & 15.808 & 164.02 & 15.168 & 0.0422 & \\ 
			    & 36.828 & 161.93 & 35.265 & 0.0443 & \\ 
			\hline
		\end{tabular}
		\caption{\sl Numerical results for the stochastic Black--Scholes equation with multiplicative noise \eqref{eq:BS}.
		}
		\label{table:black-mult}
	\end{center}
\end{table}

To obtain approximations of the function $x \mapsto X^z_T(x)$ for a given realization
$z$ of $W$, we chose $\xi^{n,m,j}$, $n \in \{1, \dots, N \}$, $m \in \{0, \dots, M-1\}$,
$j \in \N$, i.i.d.\ uniformly distributed over $[95,105]^d$ and set 
$M=25000$ together with 
$\gamma_m = 10^{-1} \mathbbm{1}_{[0,5000]}(m) 
+ 10^{-2} \mathbbm{1}_{(5000,12000]}(m)
+ 10^{-3} \mathbbm{1}_{(12000,20000]}(m)
+ 10^{-4} \mathbbm{1}_{(20000,25000)}(m)$, $m \in \{0,1,\dots, M-1\}$. Table 
\eqref{table:black-mult-domain} reports Monte Carlo estimates of the $L^1$-error
\[
\frac{1}{10^d} \int_{[95,105]^d} 
  \bigl|
    \mathbb{V}_N^{1,\mathbb{S}^{z,N}_M}\!\!(\Theta^{z,N}_M,x)-X^z_T(x)
  \bigr|
  \,
  dx,
\]
the $L^1$-norm of the relative errors
\[
  \frac{1}{10^d} 
  \int_{[95,105]^d} 
  \frac{
    \bigl|
      \mathbb{V}_N^{1,\mathbb{S}^{z,N}_M}\!\!(\Theta^{z,N}_M,x) - X^z_T(x) 
	\bigr|
  }{
    | X^z_T(x) |
  } \, dx,
\]
the $L^2$-error
\[
  \brak{
    \frac{1}{10^d} \int_{[95,105]^d} 
    {
      \bigl|
        \mathbb{V}_N^{1,\mathbb{S}^{z,N}_M}\!\!(\Theta^{z,N}_M,x)-X^z_T(x)
	  \bigr|^2
	}\,dx 
  }^{1/2},
\]
the $L^2$-norm of the relative errors
\[
  \brak{
    \frac{ 1 }{ 10^d } 
    \int_{ [ 95, 105 ]^d } 
	\brak{
	  \frac{
	    \mathbb{V}_N^{ 1, \mathbb{S}^{ z, N }_M}\!\!( \Theta^{ z, N }_M, x ) - X^z_T( x ) 
      }{ X^z_T( x ) } 
	}^{ \!\! 2 } \, dx
  }^{ 1 / 2 }
\]
and average computation times.

\begin{table}
	\begin{center}
		\small
		\begin{tabular}{|c|c|c|c|c|c|}
			\hline
			$ d $ & $L^1$-error & \makecell{$L^1$-norm of\\ rel. errors} 
			& $L^2$-error & \makecell{$L^2$-norm of\\
			rel. errors} & \makecell{Avg.\ runtime\\ in seconds} \\
			\hline
			  & 0.430716 & 0.002592 & 0.456547 & 0.002743 & \\
			  & 0.622716 & 0.003764 & 0.641510 & 0.003868 & \\
			1 & 0.216329 & 0.001250 & 0.258247 & 0.001486 & 410.3 \\
			  & 0.272603 & 0.006944 & 0.278424 & 0.007048 & \\
			  & 0.241477 & 0.005889 & 0.245629 & 0.005968 & \\
			\hline
			  & 0.022888 & 0.002038 & 0.028332 & 0.002548 & \\
			  & 0.112078 & 0.006504 & 0.115706 & 0.006765 & \\
			5 & 0.102526 & 0.004616 & 0.114451 & 0.005101 & 465.4 \\
			  & 0.122803 & 0.006738 & 0.129512 & 0.007038 & \\
			  & 0.062811 & 0.005436 & 0.067419 & 0.005876 & \\
			\hline
			   & 0.179861 & 0.003317 & 0.209160 & 0.003888 & \\
			   & 0.103214 & 0.004125 & 0.119847 & 0.004737 & \\
			10 & 0.057559 & 0.004947 & 0.061491 & 0.005252 & 514.0 \\
			   & 0.144633 & 0.003024 & 0.166896 & 0.003522 & \\
			   & 0.041878 & 0.003073 & 0.051966 & 0.003885 & \\
			\hline
			   & 0.177683 & 0.006380 & 0.188235 & 0.006810 & \\
			   & 0.086940 & 0.003106 & 0.104269 & 0.003784 & \\
			20 & 0.470434 & 0.016304 & 0.475733 & 0.016575 & 633.1 \\
			   & 0.089995 & 0.005535 & 0.094694 & 0.005873 & \\
			   & 0.159848 & 0.001588 & 0.200253 & 0.001995 & \\
			\hline
			   & 0.340198 & 0.008492 & 0.350218 & 0.008784 & \\
			   & 0.215024 & 0.017718 & 0.217900 & 0.018024 & \\
			50 & 0.553019 & 0.024297 & 0.554953 & 0.024431 & 1036.3 \\
			   & 0.174656 & 0.008763 & 0.178200 & 0.008964 & \\
			   & 0.657705 & 0.009544 & 0.667429 & 0.009705 & \\
			\hline
			    & 0.113876 & 0.002180 & 0.142040 & 0.002715 & \\
			    & 0.441801 & 0.019108 & 0.443193 & 0.019189 & \\
			100 & 0.279611 & 0.012461 & 0.283533 & 0.012669 & 1761.2 \\
			    & 0.231222 & 0.007541 & 0.237893 & 0.007782 & \\
			    & 0.313630 & 0.010046 & 0.318047 & 0.010209 & \\
			\hline
		\end{tabular}
		\caption{\sl $L^1$- and $L^2$-errors for the stochastic Black--Scholes equation with 
		multiplicative noise \eqref{eq:BS}.
		}
		\label{table:black-mult-domain}
	\end{center}
\end{table}

\subsection{Zakai equations}
\label{subsec:Zakai}
In this subsection we apply our approximation method to  Zakai equations of 
the form \eqref{eq:Zakai} below. Consider a signal process 
$S \colon [0,T] \times \Omega \to \mathbb{R}^d$ satisfying 
$\P(S_0\in A)=\int_A \varphi(x)\,dx$, $A \in \mathcal{B}(\R^d)$, for an
initial density $\varphi \colon \R^d \to \R_+$ given by 
\[
\varphi(x) = \exp(- \pi \Vert x\Vert_{\R^d}^2), \quad x \in \R^d,
\]
with dynamics
\[
S_t = S_0 + \int_0^t \kappa(S_s)\,ds +  \sigma W_t, \quad 
t \in [0,T],
\]
for $\kappa \colon \R^d \to \R^d$ of the form
\[
\kappa(x) = \frac{\beta x}{\brak{1+\Vert x\Vert_{\R^d}^2}}, \quad x \in \R^d,
\]
a $d \times d$-matrix $\sigma$ with components $\sigma_{i,j} = d^{- \nicefrac{1}{2}}$ for 
all $1 \le i,j \le d$ and a standard $d$-dimensional Brownian motion 
$W \colon [0,T] \times \Omega \to \R^d$ 
independent of $S_0$. Let the observation process $Z \colon [0,T] \times \Omega \to \R^d$
be given by 
\be \label{ZV}
Z_t = \int_0^t h(S_s) ds + V_t, \quad t \in [0,T]
\ee
for a function $h \colon \R^d \to \R^d$ of the form
$h(x) = x$, $x \in \R^d,$ and a
standard $d$-dimensional Brownian motion $V \colon [0,T] \times \Omega \to \R^d$
independent of $(S_0, W)$.
Then the Zakai equation 
 \begin{align}
&X_t(x) 
\label{eq:Zakai}\\
&= 
\varphi(x) 
+ 
\int_0^t  
\bigg(
\tfrac{1}{2} {\textstyle\sum\limits_{i,j=1}^d} \tfrac{ \partial^2}{\partial{x_i}\partial{x_j}}X_s (x)
- \textstyle{\sum\limits_{i=1}^d} \tfrac{ \partial}{\partial{x_i}}\big(\kappa_i(x) X_s(x)\big)
\bigg) ds + \int_0^t X_s(x) \, \langle h(x),  dZ_s \rangle_{\R^d} \nonumber
\end{align}
describes an unnormalized conditional density of $S_t$ along a realization of $(Z_s)_{s \in [0,t]}$, 
that is, one has 
\[
\frac{\int_{A} X_t(x) dx}{\int_{\R^d} X_t(x) dx} = \P \edg{S_t \in A \mid (Z_s)_{s \in [0,t]}}
 \quad \mbox{for every Borel subset $A \subseteq \R^d$ and all $t \in [0,T]$;}
\]
see \cite{zakai1969optimal}. Since 
 \[
 \textstyle{\sum\limits_{i=1}^d} \tfrac{ \partial}{\partial{x_i}}\big(\kappa_i(x) X_s(x)\big) =
 \textstyle{\sum\limits_{i=1}^d} X_s(x) \tfrac{ \partial}{\partial{x_i}}\kappa_i (x)
 + \big\langle \kappa( x ), \nabla X_s (x) \big\rangle_{ \R^d },
 \]
the Zakai equation \eqref{eq:Zakai} can equivalently be written as

\begin{align}
X_t(x)
&  
=\varphi(x) 
 -
 \int_0^t  \textstyle{\sum\limits_{i=1}^d} X_s(x) \tfrac{ \partial}{\partial{x_i}}\kappa_i (x) ds
 +
 \int_0^t X_s(x) \,\langle h(x), dZ_s \rangle_{\R^d}
 \nonumber
 \\
 & \quad 
 +
 \int_{ 0 }^{ t } \bigg(
 \tfrac{1}{2} {\textstyle\sum\limits_{i,j=1}^d} \tfrac{ \partial^2}{\partial{x_i}\partial{x_j}}X_s(x)
 -  \big\langle \kappa( x ), \nabla X_s ( x ) \big\rangle_{ \R^d }
 \bigg) ds,
 \nonumber
 \end{align}
which is of the form of our standard SPDE \eqref{eq:SPDE} for $\mu(x) = - \kappa(x)$, $
f(x,u,w) = \sum_{i=1}^d u \frac{\partial}{\partial{x_i}}\mu_i (x),
$
$
b(x,u,w) = u h(x)$ for all $x,w \in \R^d$, $u \in \R$ and $\sigma_{i,j} = d^{- \nicefrac{1}{2}}$, $1 \le i,j \le d$.

To apply our method to the Zakai equation \eqref{eq:Zakai}, we
used Framework~\ref{frame:adam} with $T=0.1$, $N=10$, $\varepsilon=10^{-8}$, 
$\varphi(x) = \exp(- \pi \Vert x\Vert_{\R^d}^2)$, 
\[
\mu(x) = - \kappa(x) = - \frac{x}{4 \brak{1+\Vert x\Vert_{\R^d}^2}},
\]
$\sigma_{i,j} = d^{- \nicefrac{1}{2}}$, $1 \le i,j \le d$, $h(x) = x$,
$
f(x,u,v) = \sum_{i=1}^d u \frac{\partial}{\partial{x_i}}\mu_i (x) 
$,
$
b(x,u,v) = u h(x),
$
$\delta = d$, $Z$ given by \eqref{ZV}, 
$H(t,s,x,v)= x + \mu(x)(t-s) + \sigma v$ and
\[
\mathcal{H}_n(x,u,v,w)
=u \brak{1 + 
{\textstyle{\sum\limits_{i=1}^d}} \tfrac{\partial}{\partial{x_i}}\mu_i (x) (t_{n+1}-t_n)
+ \langle h(x), w \rangle_{\R^d}}
\] 
for $x,v,w \in \R^d$, $u \in \R$ and $s,t \in [0, T]$ such that $s < t$.

We fixed a realization $z$ of $Z$ and ran our algorithm to approximate $X^z_T(0)$ five times.
Each time we set $\xi^{n,m,j}= 0 \in \R^d$ for all $n \in \{1, \dots, N \}$, $m \in \{0, \dots, M-1\}$,
$j \in \N$ and $M=35000$ together with $\gamma_m = 10^{-2} \mathbbm{1}_{[0,10000]}(m) 
+ 10^{-3} \mathbbm{1}_{(10000,20000]}(m)
+ 10^{-4} \mathbbm{1}_{(20000,30000)}(m) 
+ 10^{-5} \mathbbm{1}_{(30000,35000)}(m)$ for $m \in \{0,1,\dots, M-1\}$.
Table \ref{table:Zakai} shows results of 
$\mathbb{V}^{1,\mathbb{S}^{z,N}_M}_N(\Theta^{z,N}_M, 0)$
for the five runs, computation times, reference solutions 
$X^z_T(0)$ obtained with the Monte Carlo method\footnote{Since the 
Monte Carlo method of \cite{MR4630241} is also random, its output slightly 
varies across different runs.} developed in \cite{MR4630241}
and pathwise relative errors 
\begin{equation} \label{relerror-Zakai}
\frac{
  \bigl|
    \mathbb{V}_N^{1,\mathbb{S}^{z,N}_M}(\Theta^{z,N}_M,0) 
    - X^z_T(0)
  \bigr|
}{
  \bigl| X^z_T(0) \bigr|
}
\end{equation}
in each of the cases $d \in \{1,5,10,20,50,100\}$. 
The last column shows averages of the relative errors
\eqref{relerror-Zakai} over the five different runs.

\begin{table}
	\begin{center}
		\small
		\begin{tabular}{|c|c|c|c|c|c|}
			\hline
			$d$ & \makecell{Result of\\ the approx.\\ algorithm} & \makecell{Runtime\\ in\\ seconds} & \makecell{Reference\\ solution} & \makecell{Relative\\ pathwise\\ error} & \makecell{Avg.\\ relative\\ error}\\
			\hline
			  & 0.7638 & 311.42 & 0.7680 & 0.0054 & \\
			1 & 0.7622 & 305.34 & 0.7683 & 0.0079 & \\
			  & 0.7636 & 307.01 & 0.7678 & 0.0054 & 0.0066 \\
			  & 0.7620 & 307.02 & 0.7677 & 0.0073 &  \\
			  & 0.7630 & 304.42 & 0.7683 & 0.0069 &  \\
			\hline
			  & 0.4343 & 326.21 & 0.4419 & 0.0173 & \\
			  & 0.4365 & 333.13 & 0.4420 & 0.0124 &  \\
			5 & 0.4354 & 333.06 & 0.4422 & 0.0154 & 0.0140 \\
			  & 0.4364 & 330.98 & 0.4422 & 0.0131 &  \\
			  & 0.4373 & 331.15 & 0.4424 & 0.0116 & \\
			\hline
			   & 0.2976 & 337.32 & 0.2980 & 0.0011 & \\
			   & 0.2989 & 335.47 & 0.2990 & 0.0003 &  \\
			10 & 0.2928 & 338.50 & 0.2986 & 0.0194 & 0.0073 \\
			   & 0.2944 & 338.45 & 0.2981 & 0.0127 &  \\
			   & 0.2978 & 335.78 & 0.2986 & 0.0027 &  \\
			\hline
			   & 0.1786 & 346.75 & 0.1849 & 0.0336 & \\
			   & 0.1778 & 344.19 & 0.1852 & 0.0401 &  \\
			20 & 0.1814 & 344.19 & 0.1848 & 0.0187 &  0.0298 \\
			   & 0.1793 & 344.59 & 0.1852 & 0.0317 &  \\
			   & 0.1804 & 344.90 & 0.1850 & 0.0248 &  \\
			\hline
			   & 0.0721 & 364.12 & 0.0719 & 0.0026 &  \\
			   & 0.0719 & 367.42 & 0.0719 & 0.0002 &  \\
			50 & 0.0713 & 367.33 & 0.0721 & 0.0109 &  0.0209 \\
			   & 0.0721 & 370.27 & 0.0717 & 0.0057 &  \\
			   & 0.0715 & 367.57 & 0.0720 & 0.0064 &  \\
			\hline
			    & 0.0262 & 385.54 & 0.0273 & 0.0383 &  \\
			    & 0.0267 & 397.91 & 0.0270 & 0.0097 &  \\
			100 & 0.0266 & 399.64 & 0.0270 & 0.0175 & 0.0193 \\
			    & 0.0266 & 397.15 & 0.0271 & 0.0182 &  \\
			    & 0.0266 & 403.75 & 0.0270 & 0.0130 &  \\
			\hline
		\end{tabular}
		\caption{\sl Numerical results for the Zakai equation \eqref{eq:Zakai}.
		}
		\label{table:Zakai}
	\end{center}
\end{table}

To approximate the function $x \mapsto X^z_T(x)$, we chose $\xi^{n,m,j}$,
 $n \in \{1, \dots, N \}$, $m \in \{0, \dots, M-1\}$, $j \in \N$, i.i.d.\ uniformly distributed on $[-1,1]^d$ 
 and set $M=35000$ together with 
$\gamma_m = 10^{-2} \mathbbm{1}_{[0,10000]}(m) 
+ 10^{-3} \mathbbm{1}_{(10000,20000]}(m)
+ 10^{-4} \mathbbm{1}_{(20000,30000)}(m) + 10^{-5} \mathbbm{1}_{(30000,35000)}(m)$ for $m \in \{0, \dots, M-1\}$.
Table \ref{table:Zakai-domain} shows Monte Carlo estimates of the $L^1$- and $L^2$-errors \eqref{L1error} 
and \eqref{L2error} for a given realization $z$ of $Z$ and five runs of our approximation algorithm 
together with average computation times in each of the cases $d \in \{1,5,10,20,50,100\}$.

\begin{table}
	\begin{center}
		\small
		\begin{tabular}{|c|c|c|c|}
			\hline
			$ d $ & $L^1$-error & $L^2$-error & \makecell{Avg.\ runtime\\ in seconds}\\
			\hline
			& 0.0245668 & 0.0251878 & \\
			& 0.0201828 & 0.0206923 & \\
			1 & 0.0181104 & 0.0189806 & 370.3\\
			& 0.0187391 & 0.0236999 & \\
			& 0.0192202 & 0.0200646 & \\
			\hline
			& 0.0018356 & 0.0032205 & \\
			& 0.0027927 & 0.0053121 & \\
			5 & 0.0027895 & 0.0051387 & 463.3\\
			& 0.0025232 & 0.0047673 & \\
			& 0.0031986 & 0.0062190 & \\
			\hline
			& 0.0003742 & 0.0012530 & \\
			& 0.0002345 & 0.0009120 & \\
			10 & 0.0003458 & 0.0012521 & 555.2\\
			& 0.0002777 & 0.0010778 & \\
			& 0.0003660 & 0.0014192 & \\
			\hline
			& 0.0000024 & 0.0000108 & \\
			& 0.0000026 & 0.0000227 & \\
			20 & 0.0000016 & 0.0000172 & 738.3\\
			& 0.0000020 & 0.0000194 & \\
			& 0.0000031 & 0.0000142 & \\
			\hline
			& 0.0000027 & 0.0000027 & \\
			& 0.0000008 & 0.0000009 & \\
			50 & 0.0000004 & 0.0000005 & 1251.1\\
			& 0.0000005 & 0.0000005 & \\
			& 0.0000010 & 0.0000010 & \\
			\hline
			& 0.0000003 & 0.0000004 & \\
			& 0.0000005 & 0.0000005 & \\
			100 & 0.0000003 & 0.0000003 & 2152.4\\
			& 0.0000009 & 0.0000010 & \\
			& 0.0000011 & 0.0000011 & \\
			\hline
		\end{tabular}
		\caption{\sl $L^1$- and $L^2$-errors for the Zakai equation \eqref{eq:Zakai}.
		}
		\label{table:Zakai-domain}
	\end{center}
\end{table}

\subsection*{Acknowledgments}
This work was partially funded by the Swiss National Science Foundation
Grant 175699 ``Higher Order Numerical Approximation Methods 
for Stochastic Partial Differential Equations'' 
and the
Nanyang Assistant Professorship Grant ``Machine Learning Based Algorithms in Finance and Insurance''.
Moreover, this work has been partially funded by the European Union (ERC, MONTECARLO, 101045811).
The views and the opinions expressed in this work are however those of the authors only and do
not necessarily reflect those of the European Union or the European Research Council (ERC).
Neither the European Union nor the granting authority can be held responsible for them.
We also gratefully acknowledge the 
Cluster of Excellence EXC 2044-390685587, Mathematics Münster: Dynamics-Geometry-Structure 
funded by the Deutsche Forschungsgemeinschaft (DFG, German Research Foundation).
%
%
\bibliographystyle{acm}

\end{document}